\documentclass[11pt]{article}
\usepackage{epsfig,latexsym,graphicx,amssymb}
\usepackage[dvips]{psfrag}

\usepackage{amssymb, amsmath, amscd, amsthm, amstext, amsfonts}

\def\AA{{\mathbb A}} \def\BB{{\mathbb B}} \def\CC{{\mathbb C}}
  \def\DD{{\mathbb D}}

 \def\NN{{\mathbb N}} 

 \def\RR{{\mathbb R}} 

\def\UU{{\mathbb U}} \def\VV{{\mathbb V}} 

 \def\ZZ{{\mathbb Z}}

\def\mfT{{\mathfrak{T}}}

\def\Si{\Sigma}
\def\La{\Lambda}
\def\De{\Delta}
\def\Om{\Omega}
\def\Ga{\Gamma}

\def\la{\lambda}
\def\be{\beta}
\def\de{\delta}

\def\vfi{\varphi}
\def\ve{\varepsilon}

\def\diffM{\mbox{{\rm Diff\,}}^1(M)}
\def\diffrM{\mbox{{\rm Diff\,}}^r(M)}
\def\Perf{\mbox{{\rm Per}}(f)}

\def\Perfn{\mbox{{\rm Per}}_n (f)}

\def\cA{{\cal A}}  \def\cG{{\cal G}}  \def\cS{{\cal
S}}
\def\cB{{\cal B}}   \def\cN{{\cal N}} \def\cT{{\cal
T}}
\def\cC{{\cal C}}   \def\cO{{\cal O}} \def\cU{{\cal
U}}
\def\cD{{\cal D}}    \def\cV{{\cal
V}}
    \def\cW{{\cal
W}}
   \def\cR{{\cal R}}

\newcommand{\eqdef}{\stackrel{\scriptscriptstyle\rm def}{=}}

\newtheorem{theo}{Theorem}
\newtheorem{coro}{Corollary}

\newtheorem{lemm}{Lemma}[section]

\newtheorem{defi}[lemm]{Definition}
\newtheorem{prop}[lemm]{Proposition}
\newtheorem{theor}[lemm]{Theorem}
\newtheorem{clai}[lemm]{Claim}
\newtheorem{conj}{Conjecture}

\newtheorem{ques}{Question}

\newtheorem{rema}[lemm]{Remark}

\newenvironment{demo}[1][Proof]{{\medskip
\noindent\bf #1:
}}{\hfill$\Box$\medskip}

\def\st{{\operatorname{s}}}
\def\sst{{\operatorname{ss}}}
\def\cst{{\operatorname{cs}}}

\def\ct{{\operatorname{c}}}
\def\cut{{\operatorname{cu}}}
\def\ut{{\operatorname{u}}}
\def\uut{{\operatorname{uu}}}

\def\inds{{\mbox{\rm{ind}}^{\,\st}}}
\def\indu{{\mbox{\rm{ind}}^{\,\ut}}}
\def\loc{{\operatorname{loc}}}

\title{Abundance of $C^1$-robust homoclinic tangencies}
\author{Christian Bonatti and  Lorenzo J. D\'\i az\thanks{This paper
was partially supported by CNPq, Faperj, and PRONEX (Brazil) and
the Agreement in Mathematics Brazil-France. We acknowledge the
warm hospitality of I.M.P.A, Institute de Math\'ematiques de
Bourgogne, and PUC-Rio during the stays while preparing this
paper}}

\begin{document}
\maketitle

\bigskip

\begin{flushright}{{\tt{To Carlos Gutierrez (1944 -- 2008), in memoriam}}}
\end{flushright}

\bigskip

\begin{abstract}
A diffeomorphism $f$ has a $C^1$-robust homoclinic tangency if there
is a $C^1$-neighbourhood $\cU$ of $f$ such that every
diffeomorphism in $g\in \cU$ has a hyperbolic set $\La_g$, depending
continuously on $g$, such that the stable and unstable manifolds of
$\La_g$ have some non-transverse intersection. For every manifold of
dimension greater than or equal to three, we exhibit a local
mechanism (blender-horseshoes) generating diffeomorphisms with
$C^1$-robust homoclinic tangencies.

Using blender-horseshoes, we prove that  homoclinic classes of
$C^1$-generic diffeomorphisms containing saddles with different
indices and that do not admit  dominated splittings (of appropriate
dimensions) display $C^1$-robust homoclinic tangencies.

\medskip

\noindent{\bf keywords:}
chain recurrence set, dominated splitting,
heterodimensional cycle,
homoclinic class, homoclinic tangency,
hyperbolic set.

\medskip

\end{abstract}

\medskip

\section{Introduction}

\subsection{Framework and general setting}
A homoclinic tangency is a dynamical mechanism which is  at the
heart of a great variety  of non-hyperbolic phenomena: persistent
coexistence of infinitely many sinks \cite{NeIHES}, H\'enon-like
strange attractors \cite{BeCaAM,MV}, super-exponential growth of
the number of periodic points \cite{KaAnMa}, and non-existence of
symbolic extensions \cite{DoNeIn}, among others. Moreover,
homoclinic bifurcations (homoclinic tangencies and
heterodimensional cycles) are conjectured to be the main source of
non-hyperbolic dynamics (Palis denseness conjecture, see
\cite{PaAst}).

In this paper, we present a local mechanism generating $C^1$-robust
homoclinic tangencies. Using this construction, we show that the
occurrence of robust tangencies is a quite general phenomenon in the
non-hyperbolic setting, specially when the dynamics does not admit a
suitable dominated splitting.

Let us now give some basic definitions (in Section~\ref{s.definitions}, we will
state precisely the definitions involved in
this paper).
A  transitive hyperbolic set $\La$ has a \emph{homoclinic tangency}
if there is a pair of points  $x,y\in \La$ such that the stable leaf
$W^\st(x)$ of $x$ and the unstable leaf $W^\ut(y)$ of $y$  have some
non-transverse intersection

Given a hyperbolic set $\La$ of a diffeomorphism $f$, for $g$
close to $f$, we denote by $\La_g$ the hyperbolic set of $g$ which
is the {\emph{continuation}} of $\La$ (i.e., $\La_g$ is close to
$\La$ and the dynamics of $f$ on $\La$ and $g$ on $\La_g$ are
conjugate).

\begin{defi}[Robust cycles] $\,$
\begin{itemize}
\item
{\em Robust homoclinic tangencies:} A transitive hyperbolic set
$\La$ of a $C^r$-diffeo\-mor\-phism $f$ has  a {\emph{$C^r$-robust
homoclinic tangency}} if there is a $C^r$-neighborhood $\cN$ of
$f$ such that for every $g\in \cN$ the continuation $\La_g$ of
$\La$ for $g$ has a homoclinic tangency.
\item
{\em{Robust heterodimensional cycles:}} A diffeomorphism $f$ has a
\emph{$C^r$-robust heterodimensional cycle} if there are
transitive hyperbolic sets $\La$ and $\Sigma$ of $f$ whose stable
bundles have different dimensions and a  $C^r$-neighborhood $\cV$
of $f$
  such that $W^\st(\La_g) \cap W^\ut (\Sigma_g)\ne \emptyset$ and
$W^\ut(\La_g) \cap W^\st (\Sigma_g)\ne \emptyset$, for every
diffeomorphism $g\in \cV$.
\end{itemize}
\label{d.robusthomoclinic}
\end{defi}

Note that, by Kupka-Smale theorem, $C^r$-generically, invariant
manifolds of periodic points are in general position. Hence,
generically, the non-transverse intersections in a robust cycle
(tangency or heterodimensional cycle) involve non-periodic points
(i.e., at least a non-trivial hyperbolic set).

\medskip

In \cite{NeProvidence}, Newhouse constructed surface
diffeomorphisms having hyperbolic sets (called {\emph{thick
horseshoes}}) exhibiting $C^2$-robust homoclinic tangencies.
Later, he proved that, in dimension two, homoclinic tangencies of
$C^2$-diffeomorphisms  yield thick horseshoes with $C^2$-robust
homoclinic tangencies, \cite{NeIHES} (see also \cite{PaTabook} for a
broad discussion of homoclinic bifurcations on surfaces). With the
same $C^2$-regularity assumption, theorems in \cite{RoETDS,PaViAnMa}
extend Newhouse result, proving that  homoclinic tangencies in
any dimension lead to $C^2$-robust homoclinic tangencies.
In this paper, we study the occurrence of robust homoclinic tangencies in the $C^1$-setting.

Newhouse construction (thick horseshoes with robust tangencies)
involves distortion estimates which are typically $C^2$. The
results in \cite{UrAENS} present some obstacles for carrying this
construction to the $C^1$-topology: $C^1$-generic surface
diffeomorphisms do not have thick horseshoes. Recent results by
Moreira in \cite{Gugu} are a strong indication that there are no
surface diffeomorphisms exhibiting $C^1$-robust homoclinic
tangencies\footnote{This question is closely related to the open
problem of $C^1$-density of hyperbolic diffeomorphisms on compact
surfaces (Smale's density conjecture). In fact, Moreira's result
imply that there are no $C^1$-robust homoclinic tangencies
associated to hyperbolic {\emph{basic}} sets of surface
diffeomorphisms. See \cite{ABCDFund} for a discussion of the
current state of this conjecture.}.

Nevertheless, in higher dimensions, there are examples of
diffeomorphisms having hyperbolic sets with $C^1$-robust
tangencies. For instance, the product of a non-trivial hyperbolic
attractor by a normal expansion gives a hyperbolic set $\La$ of
saddle type, whose stable manifold has a  topological dimension
greater than the dimension of its stable bundle. Then the set
$\La$ can play the role of thick horseshoes in Newhouse
construction. Geometrical constructions using these kind of
``thick" hyperbolic sets provide examples of systems with
\emph{$C^1$-robust heterodimensional cycles}\footnote{A
{\emph{heterodimensional cycle}} is a cycle associated to saddles
having different indices. } (see \cite{AbSm}) or $C^1$-robust
tangencies (see \cite{Simon,AsPAMS}). But these constructions
involve quite specific global dynamical configurations, thus they
cannot translate to a general setting.

\subsection{Robust homoclinic tangencies}

The aim of this paper is to show that the existence of
$C^1$-robust homoclinic tangencies is a common phenomenon in the
non-hyperbolic setting. For instance, next result is a consequence
of the local mechanism for robust tangencies in
Theorem~\ref{t.tangency}.

\begin{theo}
\label{t.main} Let $M$ be a compact manifold with $\dim (M)\geq
3$. There is a residual subset $\cR$ of $\diffM$ such that, for
every $f\in\cR$ and every periodic saddle $P$ of $f$ such that
\begin{itemize}
\item {\bf (index variability)} the homoclinic class $H(P,f)$  of $P$ has a periodic saddle $Q$
with  $\dim (E^\st(Q))\neq \dim (E^\st(P))$,
\item {\bf (non-domination)}
the stable/unstable splitting $E^\st(R)\oplus E^\ut(R)$ over the set of 
saddles $R$ homoclinically related with $P$ is not dominated,
\end{itemize}
then the saddle $P$ belongs to a transitive hyperbolic set having a
$C^1$-robust homoclinic tangency.
\end{theo}

For the precise definitions of homoclinic class and dominated splitting
see Definitions~\ref{d.homoclinic} and \ref{d.dominate}.
Let us reformulate  Theorem~\ref{t.main} by focusing on the homoclinic 
class  of a prescribed periodic orbit:

\begin{coro}\label{c.main}
Let $M$ be a compact manifold with $\dim (M)\geq
3$. Consider a diffeomorphism $f$ with a saddle $P_f$ whose
continuation $P_g$ is defined for all $g$ in a neighborhood $\cU$ of $f$ in
$\diffM$. Assume that
\begin{itemize}
\item {\bf (generic index variability)}
there is a  residual subset $\cG$ of $\cU$ such that, for
every $g\in\cG$, the homoclinic class of  $P_g$ of $f$ contains a saddle $Q$ of different index,
\item {\bf (robust non-domination)} 
for every $g\in \cU$, the stable/unstable splitting $E^\st(R)\oplus E^\ut(R)$ over the set of 
saddles $R$ homoclinically related with $P_g$ is not dominated.
\end{itemize}
Then there is an open and dense subset $\cC$ of $\cU$ of diffeomorphisms
$g$ such that
the saddle $P_g$ belongs to a transitive hyperbolic set with a
$C^1$-robust homoclinic tangency. 
\end{coro}

\begin{rema}
 The diffeomorphisms $f$ in the residual subset $\cR$ of $\diffM$
in Theorem~\ref{t.main}
satisfy the following properties (see
{\em{\cite[Section 2.1]{ABCDW}}} and
{\em{\cite[Appendix B.1.1]{BDVbook}}}):
\begin{itemize}
\item
Every homoclinic class $H(P_f,f)$ of $f$ depends continuously on
$f\in \cR$. Therefore, if $H(P_f,f)$ has a dominated splitting
then $H(P_g,g)$ also has a dominated splitting whose bundles have
constant dimension for all $g\in \cR$ close to $f$.

\item
Assume that a
   homoclinic class  $H(P_f,f)$ of $f\in \cR$ contains
saddles of stable indices $j$ and $k$, $j\ne k$. Then the
homoclinic class $H(P_g,g)$ also contains saddles of stable
indices $j$ and $k$ for every $g\in \cR$ close to $f$.
\end{itemize}
In other words, the conditions in Theorem~\ref{t.main} are
$C^1$-open in the residual set $\cR$ of $\diffM$.
\end{rema}

The  \emph{index interval}  of a homoclinic class $H$ is the
interval $[i,j]$, where $i$ and $j$ are the minimum and the
maximum of the \emph{$\st$-indices} (dimension of the stable
bundle) of the periodic points in $H$. The homoclinic
class $H$ has \emph{index variation} if $i<j$. Given a transitive hyperbolic set $\La$ its
\emph{$\st$-index} is the dimension of its stable bundle.

\begin{coro} For every
 diffeomorphism $f$  in the residual subset $\cR$
of $\diffM$, any homoclinic class $H$ of $f$ with index variation,
and every $k\in [i,j]$, where $[i,j]$ is the index interval of $H$, one
has:
\begin{itemize}
\item either there is a dominated splitting $E\oplus_{_<} F$ (i.e., $F$ dominates $E$) with
$\dim (E)=k$,
\item or there is a hyperbolic transitive set $\La\subset H$ with
$\st$-index $k$ having a
$C^1$-robust homoclinic tangency.
\end{itemize}
\label{c.index}
\end{coro}

When we are interested only in the existence of robust homoclinic tangencies, without paying attention to the index
of the hyperbolic set involved in their generation, there is the following reformulation:

\begin{coro}
There is a residual subset $\cG$ of $\,\diffM$ 
such that for every diffeomorphism $f\in \cG$ and every homoclinic class  
$H(P,f)$ of $f$ with index  interval $[i,j]$, $j>i$,

\begin{itemize}
 \item 
either $H(P,f)$ has a dominated splitting 
$$
T_{H(P,f)}M= E^\cst\oplus_{_<}
E_1\oplus_{_<}\cdots\oplus_{_<}E_j\oplus_{_<}E^\cut,
$$ 
where
$\dim(E^\cst)=i$ and
$E_1,\dots,E_j$ are one-dimensional,
\item
or the homoclinic class $H(P,f)$ contains a 
transitive
hyperbolic set with a robust homoclinic tangency.
\end{itemize}
\label{c.dichotomy}
\end{coro}

In the first case of Corollary~\ref{c.dichotomy}, 
we say that $H(P,f)$ has an {\emph{indices adapted dominated splitting.}}

\medskip

The previous results have an interesting formulation  for  {\emph{tame
diffeomorphisms,}} i.e., the $C^1$-open set $\cT(M)$ of $\diffM$ of
diffeomorphisms having finitely many chain recurrence classes (see
Definition~\ref{d.chain}) in a robust way. 
We define $\cW(M)\eqdef\diffM \setminus \overline{\cT(M)}$ as the set
of {\emph{wild diffeomorphisms.}}
Let us observe
that, for an open and dense subset of $\cT(M)$, a chain recurrence
class is either hyperbolic or has index variation, see
\cite{ABCDW}. 

Given a chain recurrence class $C$ of $f$ we first consider the
{\emph{finest dominated splitting}} over $C$ (i.e., the bundles of
this splitting can not be decomposed in a dominated way). Then we
let $E^\st$ (resp. $E^\ut$) be the sum of the uniformly
contracting (resp. expanding) bundles of this splitting (these bundles may be 
trivial,
see \cite{BV}). The
bundles $E_1,\dots, E_k$ are the remaining non-hyperbolic bundles
of the finest dominated splitting of $C$. In this way, we get a dominated
splitting over $C$
$$
T_CM= E^\st\oplus_{_<}
E_1\oplus_{_<}\cdots\oplus_{_<}E_k\oplus_{_<}E^\ut,
$$ where
$E^\st$ and $E^\ut$ are uniformly contracting and expanding, and
$E_1,\dots,E_k$ are indecomposable and non-hyperbolic. We call
this splitting the {\emph{finest central dominated splitting}} of
the chain recurrence $C$.

\begin{rema}
Let $f$ be any tame
diffeomorphisms and $H(P,f)$ any homoclinic class of $f$  which is 
far from robust homoclinic tangencies. Then the 
finest central dominated splitting of $H(P,f)$ 
is indices adapted. For tame diffeomorphisms, the corollary below gives a more precise description of
the relation between the finest central dominated splitting and the robust homoclinic
tangencies associated to a homoclinic class. 
\end{rema}

\begin{coro}\label{c.tame} There is a $C^1$-open and dense subset $\cO$ of
the set $\cT(M)$ of tame diffeomorphisms such that, for every $f$ in
$\cO$ and every chain recurrence class $C$ of $f$ whose finest
central dominated splitting is
$$
T_CM= E^\st\oplus_{_<}
E_1\oplus_{_<}\cdots\oplus_{_<}E_k\oplus_{_<}E^\ut,
$$
then, for every $i=1,\dots,k$,
$$
\dim (E_i)>1 \Longleftrightarrow \left\{
\begin{array}{ll}

& \mbox{for all $j\in\{1,\dots, \dim(E_i)-1\}$,}\\
& \mbox{there is a transitive hyperbolic set $K$ of $\st$-index}\\
&\inds (K)=\dim\left(E^\st\oplus
E_1\oplus \cdots \oplus E_{i-1}\right)+j\\
&\mbox{having a $C^1$-robust homoclinic tangency.}
\end{array}
 \right.
$$
\end{coro}

\begin{rema}$\,$
\begin{enumerate}
\item
Under the hypotheses  of Corollary~\ref{c.tame}, {\em
\cite[Theorem 1.14]{BDijmj}} implies that (choosing appropriately
the open and dense subset $\cO$ of $\cT(M)$) the hyperbolic set $K$
with a $C^1$-robust homoclinic tangency is also involved in a
$C^1$-robust heterodimensional cycle.
\item
Corollary~\ref{c.tame} can also be stated for isolated chain
recurrence classes of $C^1$-generic diffeomorphisms \footnote{By
$C^1$-generic diffeomorphisms we mean diffeomorphisms in a
residual subset of $\diffM$.}.
\end{enumerate} \label{r.main}
\end{rema}

This article proceed a program for studying the generation of
{\emph{robust cycles}} (homoclinic tangencies and heterodimensional cycles)
in the $C^1$-topology. In  \cite{BDijmj} we
proved that homoclinic
classes containing periodic points with different indices
generate (by arbitrarily small $C^1$-perturbations)
$C^1$-robust heterodimensional cycles. Here
 we show that these robust heterodimensional cycles generate
 {\emph{blender-horseshoes,}} a sort of hyperbolic basic sets with
 geometrical properties resembling the thick horseshoes, see Section~\ref{ss.blenderhorseshoe} and Theorem~\ref{th.genericblender}.
We next see that, in the context of critical dynamics (some
suitable non-domination property), blender-horseshoes yield
$C^1$-robust tangencies, see Theorem~\ref{t.tangency}. In fact,
the definition and construction of blender-horseshoes (a special
class of $\cut$-blenders defined in \cite{BoDiAnMa})
 and
Theorem~\ref{t.tangency} are the technical heart of our arguments
and the main novelty of this paper.

\medskip

The results in this paper and the ones in
\cite{BDijmj} support the following conjecture:

\begin{conj}[Bonatti, \cite{umalca}]
Every $C^1$-diffeomorphism can be $C^1$-approximated either by a
hyperbolic diffeomorphism (Axiom A and no-cycle property) or by a
diffeomorphism exhibiting a $C^1$-robust cycle (homoclinic
tangency or heterodimensional cycle). \label{c.bonatti}
\end{conj}

This conjecture is a stronger version of the denseness conjecture
by Palis in \cite{PaAst} (dichotomy hyperbolicity versus
approximation by diffeomorphisms with homoclinic bifurcations).
The novelty here is that the conjecture considers two disjoint
open sets whose union is dense in the whole set of
$C^1$-diffeomorphisms: the hyperbolic ones and those with robust
cycles. In the setting of tame diffeomorphisms,  a strong version of Conjecture~\ref{c.bonatti} was proved
in \cite[Theorem 1.14]{BDijmj}: every
tame diffeomorphism can be $C^1$-approximated either by hyperbolic
diffeomorphisms or by diffeomorphisms exhibiting  robust heterodimensional
cycles. Recall that Palis conjecture for surface
$C^1$-diffeomorphisms was proved in \cite{PuSaAnMa} (due to
dimension deficiency, for surface diffeomorphisms the conjecture
only involves homoclinic tangencies).

\subsection{Newhouse domains}
 Following \cite{KaAnMa}, we
say that an open set $\cN$ of $\diffrM$ is a {\emph{$C^r$-Newhouse
domain}} if there is a dense subset $\cD$ of $\cN$ such that every
$g\in \cD$ has a homoclinic tangency (associated to some saddle).
A preliminary step  toward Conjecture~\ref{c.bonatti}  is
the following question.

\begin{ques}
Let $M$ be a closed manifold and $\cN$ be a $C^1$-Newhouse domain
of $\diffM$. Are the diffeomorphisms having $C^1$-robust
homoclinic tangencies dense in $\cN$?
 \label{q.newhouse}
\end{ques}

If it is not possible to answer positively this question in its
full generality, it would be interesting to provide sufficient
conditions for a $C^1$-Newhouse domain to contain an open and
dense subset of diffeomorphisms with $C^1$-robust homoclinic
tangencies. If the dimension of the ambient manifold is at least
three, one may also ask about the interplay between robust homoclinic
tangencies and robust heterodimensional cycles.

We now discuss briefly Question~\ref{q.newhouse}. Before going to our setting, lets us
review the discussion in \cite{ABCDFund} about this question for
$C^1$-surface diffeomorphisms.  Let $\textrm{Hyp}^1(M)$ denote the subset of $\diffM$
consisting of Axiom A diffeomorphisms. By \cite{PuSaAnMa}, for
surface diffeomorphisms, the open set
$$
\cN^1(M^2) \eqdef \textrm{Diff}^1(M^2)\setminus
\overline{\textrm{Hyp}^1(M^2)}
$$
is a Newhouse domain. The set $\cN^1(M^2)$ is the union of the closure of three
pairwise disjoint
 open
sets $\cO_1(M^2),\cO_2(M^2)$, and $\cO_3(M^2)$ defined as follows.

\begin{itemize}
\item
The set $\cO_1(M^2)$ consists of diffeomorphisms having
$C^1$-robust homoclinic tangencies.
\item
There is a residual subset $\cR_2(M^2)$ of $\cO_2(M^2)$ such that
every $f\in \cR_2(M^2)$ has a homoclinic class $H(P,f)$ that robustly does
not admit any dominated splitting. However, for every hyperbolic set
$\La$ contained in $H(P,f)$ the invariant manifolds of $\La$ meet
transversely. In this case, we say that the
diffeomorphism $f$ has a {\emph{persistently fragile homoclinic tangency}} associated
to $P$.
\item
There is a residual subset $\cR_3(M^2)$ of $\cO_3(M^2)$ such that
for every diffeomorphism $f\in \cR_3(M^2)$ and every (hyperbolic)
periodic point $P$ of $f$ the homoclinic class $H(P,f)$ is
hyperbolic. But there is a sequence of periodic points $(P_n)_n$
of $f$ such that the hyperbolic homoclinic classes $H(P_n,f)$
accumulate (Hausdorff limit) to an aperiodic class (i.e., a
recurrence class without periodic points).
\end{itemize}

As mentions above, Moreira's result in \cite{Gugu}
provides strong evidences suggesting that 
$\cO_1(M^2)$ is empty. On the other hand,
we do not know if the sets  $\cO_2(M^2)$
and $\cO_3(M^2)$ are empty or not. In fact, Smale density
conjecture (hyperbolic diffeomorphisms are dense in
$\textrm{Diff}^1(M^2)$) is equivalent to prove that these three
sets are empty.

We now explain how  the discussion above is translated to higher
dimensions. As before, we first consider non-hyperbolic diffeomorphisms, that is,
the set $\textrm{Diff}^1(M)\setminus
\overline{\textrm{Hyp}^1(M)}$.
If $\dim (M)\ge 3$ this set is not a Newhouse
domain: it contains open sets of diffeomorphisms without
homoclinic tangencies. 
Thus we consider the sets  $\textrm{Tang}^1(M)$ of diffeomorphisms having a homoclinic tangency
associated to a saddle and  $\cO_0(M)$ of 
{\emph{non-hyperbolic diffeomorphisms far from homoclinic tangencies,}}
$$
\cO_0(M)\eqdef  \diffM\setminus \overline{(\textrm{Tang}^1(M) \cup \textrm{Hyp}^1(M))}.
$$
Note that this set is not-empty and it is an open question whether it is contained
in the set of tame diffeomorphisms (in fact, the first author conjectured that $\cO_0(M)$ consists of 
tame diffeomorphisms, \cite{umalca}). The diffeomorphisms in $\cO_0(M)$ were studied
in several papers, let us just refer to \cite{Wen1,Wen2,Yang}.

From now on, we will focus on the set 
$$
\cN^1(M)\eqdef \diffM \setminus \overline{ \left( \cO_0(M)
\cup \textrm{Hyp}^1(M)\right)}.
$$ 
By definition, 
this set is a Newhouse domain. As in the case of surface diffeomorphisms, 
we split the set
$\cN^1(M)$ into three closed sets with pairwise disjoint interiors.
We first  define the set $\cO_1(M)$ similarly as  the set
$\cO_1(M^2)$, 
$$
\cO_1(M)\eqdef \{\mbox{$f\in \diffM$ with a transitive hyperbolic set with a robust 
homoclinic tangency}\}.
$$
The results in this paper implies that $\cO_1(M)$ is
non-empty, see also \cite{AsPAMS,Simon}. 

We define the set $\cO_2(M)$
by
$$
\cO_2(M)\eqdef \{f\in (\diffM \setminus \overline{\cO_1(M)}) \mbox{\,with a persistently fragile 
homoclinic tangency}\}.
$$

Consider the residual set $\cG$ of $\diffM$ in Corollary~\ref{c.dichotomy}.
Then
if $f$ is a 
diffeomorphism  in $\cG\cap \cO_2(M)$ with a  persistently fragile 
homoclinic tangency associated to $P$ then the homoclinic class
$H(P,f)$ has no index variation (otherwise one gets robust homoclinic tangencies). 

Finally, define $\cO_3(M)$ by
$$
\cO_3(M)\eqdef
\left(\diffM \setminus \overline{ \textrm{Hyp}^1(M) \cup
\cO_0(M)
\cup
\cO_1(M) \cup \cO_2(M)} \right).
$$
Corollary~\ref{c.dichotomy} implies that if 
$
f\in \cG \cap \cO_3(M)$
then every homoclinic class of $f$ has an indices adapted dominated splitting. 
The description of the accumulation of homoclinic 
classes of diffeomorphisms in $\cO_3(M)$ is a subtle issue. 
For instance, by shrinking $\cG$, for diffeomorphisms  $f\in \cG \cap \cO_3(M)$, there are $k$
and a sequence of saddles $P_n$ of index $k$ such that every $H(P_n,f)$ 
has a dominated splitting $E\oplus_{{_<}}F$ with $\dim (E)=k$ 
and the sequence of homoclinic classes $H(P_n,f)$
accumulates 
to a set $\Lambda$ that does not admit a dominated splitting  $E\oplus_{{_<}}F$ 
with $\dim(E)=k$ (the set $\Lambda$ is the Hausdorff limit of the sequence $(H(P_n,f))$).

We observe that,
as a consequence of Corollary~\ref{c.tame},
there is an open and dense subset of  $\cO_2(M)\cup \cO_3(M)$ consisting of wild diffeomorphisms.
Note that we do not known if the sets $\cO_2(M)$ and $\cO_3(M)$ are empty or
not.

Summarizing, as in the case of surface diffeomorphisms, we have that the Newhouse domain
$\cN^1(M)$ is the closure of the union of the pairwise disjoint open sets $\cO_1(M),\cO_2(M)$, and
$\cO_3(M)$.

%\begin{ques}
%Let $M$ be a closed manifold with $\dim(M)\ge 3$. State sufficient conditions
%for a
%$C^1$-Newhouse domain of $\diffM$ containing an open and dense subset
%of diffeomorphisms with $C^1$-robust heterodimensional cycles.
% \label{q.newhouseheterodimensional}
%\end{ques}

\medskip

This paper is organized as follows. In Section~\ref{s.definitions},
we recall some definitions and state some notations we will use
throughout the paper. In Section~\ref{s.bh}, we review the notion
of   $\cut$-blender in \cite{BoDiAnMa} and present the notion of
 blender-horseshoe, a key ingredient of our constructions. In
Section~\ref{s.robusttangencies}, we introduce a class of
sub-manifolds, called {\emph{folding manifolds}} relative to a
blender-horseshoe $\La$. The main result is that folding manifolds
and the local stable manifold of the blender-horseshoe $\La$ have
$C^1$-robust tangencies, see Theorem~\ref{t.tangency}. Using this
result, we state a sufficient condition for the generation of
robust homoclinic tangencies by
  homoclinic tangencies associated to
 hyperbolic sets. In
Section~\ref{s.examples},
we see that strong homoclinic intersections of
non-hyperbolic periodic points (i.e., intersections between the strong stable and unstable
manifolds) generate blender-horseshoes. We also see that such strong intersections
naturally occur in the non-hyperbolic setting. % an are associated to the unfolding of heterodimensional cycles.
Finally, in Section~\ref{s.tmain}
we  conclude the proof of  Theorem~\ref{t.main}.
We also state  a result
about the occurrence of robust heterodimensional cycles inside
non-hyperbolic chain recurrence classes, see Theorem~\ref{t.heterocycle}, which is  an extension of \cite[Theorem 1.16]{BDijmj}.

\section{Definitions and notations}
\label{s.definitions}

In this section, we define precisely the notions involved in this
paper and state some notations.

Given a closed manifold $M$, we denote by $\diffM$ the space of
$C^1$-diffeomorphisms endowed with the usual uniform topology.

A diffeomorphism $f$ has a {\emph{homoclinic tangency}} associated
to a (hyperbolic) saddle $R$ if the unstable manifold $W^\ut(R,f)$
and the stable manifold $W^\st(R,f)$ of the orbit of $R$ have some
non-transverse intersection.

The $\st$-index (resp. $\ut$-index) of a hyperbolic periodic point
$R$, denoted by $\inds(R)$ (resp. $\indu(R)$),   is the dimension
of the stable bundle $E^\st$ (resp. dimension of $E^\ut$) of
$R$. We similarly define the $\st$-index and $\ut$-index of a
transitive hyperbolic set $\La$, denoted by $\inds(\La)$ and
$\indu(\La)$, respectively.

A {\emph{heterodimensional cycle}} of a diffeomorphism $f$
consists of two hyperbolic saddles $P$ and $Q$ of $f$ of different
$\st$-indices and two heteroclinic points $X\in W^\ut(P,f)\cap W^\st(Q,f)$ and
$Y\in W^\st(P,f)\cap W^\ut(Q,f)$. In this case, we say that the
cycle is associated to $P$ and $Q$. Note that (due to insufficient
dimensions) at least one of these intersections is not transverse.
The heterodimensional cycle has {\emph{co-index $k$} if
$|\inds(Q)-\inds(P)|=k$ (note that $k\ge 1$).

\begin{defi}[Homoclinic class]
Consider a diffeomorphism $f$ and a saddle $P$ of $f$. The
{\emph{homoclinic class of $P$,}} denoted by $H(P,f)$, is the
closure of the transverse intersections of the stable and unstable
manifolds of the orbit of $P$.
\label{d.homoclinic}
\end{defi}

\begin{rema}
The homoclinic class $H(P,f)$ can be alternatively defined as the
closure of the saddles $Q$ {\emph{homoclinically related with
$P$}}: the stable manifold of the orbit of $Q$ transversely meets
the unstable manifold of the orbit of $P$ and vice-versa. Although
all  saddles homoclinically related with $P$ have the same
$\st$-index as $P$, the homoclinic class  $H(P,f)$  may
contain periodic orbits of different $\st$-index as the one of $P$
(i.e., there are homoclinic classes having index variation).
Finally, a homoclinic class is a transitive set with dense
periodic points.
\end{rema}

\begin{defi}[Chain recurrence class]
A point $x$ is \emph{chain recurrent} if for every $\varepsilon>0$
there are $\varepsilon$-pseudo orbits starting and ending at $x$.
The \emph{chain recurrence class} of $x$ for $f$, denoted by
$C(x,f)$, is the set of points $y$ such that, for every $\ve>0$,
there are $\varepsilon$-pseudo orbits starting at $x$, passing
$\ve$-close to $y$  and ending at $x$. \label{d.chain}
\end{defi}
According to \cite{BoCroIn}, for $C^1$-generic diffeomorphisms,
the chain recurrence class of any periodic point is  its
homoclinic class.

\begin{defi}[Dominated splitting]\label{d.dominate}
Consider a diffeomorphism $f$ and a compact $f$-invariant set
$\La$. A $Df$-invariant splitting $T_{\La}M=E\oplus F$ over $\La$
is {\emph{dominated}} if the fibers $E_x$ and $F_x$ of $E$ and $F$
have constant dimension and there exists $k\in \NN$ such that
$$
\frac{||D_x f^k(u)||}{||D_xf^{k} (w) ||} <  \frac{1}{2},
$$
 for
every $x\in \La$ and every pair of unitary vectors $u\in E_x$ and
$w\in F_x$.

This definition means that  vectors in the bundle $F$ are uniformly more
expanded than vectors in $E$ by the derivative $Df^k$. If it
occurs, we say that $F$ {\em dominates\/} $E$ and write
$E\oplus_{_<}F$.
\end{defi}

\begin{rema}
\label{r.severalbundles} In some cases, one needs to consider
splittings with more than two bundles. A $Df$-invariant splitting
$E_1\oplus E_2\oplus \cdots \oplus E_k$ over a set $\La$ is
dominated if for all $j\in \{1,\dots,k-1\}$ the splitting
$E_1^j\oplus E_{j+1}^k$ is dominated, where $E_i^r=E_i\oplus \cdots
\oplus E_r$, $i<r$.

We use the notation $E_1\oplus_{_<} E_2\oplus_{_<} \cdots
\oplus_{_<} E_k$, meaning that $E_{i+1}$ dominates $E_i$, or
equivalently that $E_1^j\oplus_{_<} E_{j+1}^k$.
\end{rema}

As mentioned before, the main goal of this paper is to
construct hyperbolic sets exhibiting homoclinic tangencies in a
robust way. We need the following definition.

\begin{defi}[Robust tangency]
Given a diffeomorphism $f\colon M\to
M$, a
 hyperbolic set
$\Ga$ of $f$ with a hyperbolic splitting $E^\st\oplus E^\ut$, and
a submanifold $N\subset M$ with dimension $\dim(N)=\dim(E^\ut)$,
we say that the stable manifold $W^\st(\Ga)$  of $\Ga$ and the
submanifold $N$ have a {\emph{$C^1$-robust tangency}} if for every
diffeomorphism $g$ $C^1$ close to $f$  and every submanifold $N_g$ $C^1$-close to $N$, the stable manifold
$W^\st(\Ga_g)$ of $\Ga_g$ has some non-transverse intersection
with $N_g$.
\label{d.rtm}
\end{defi}
% here $\Ga_g$ is the hyperbolic continuation of
%$\Ga$. 

We are specially interested in the case where $N$ is the
unstable manifold $W^\ut(P)$ of a periodic point $P$ of a
non-trivial hyperbolic set $\Ga$ and $N_g=W^\ut (P_g)$. In that case, one gets
$C^1$-robust homoclinic tangencies (associated to $\Ga$), recall Definition~\ref{d.robusthomoclinic}.

\bigskip

\noindent{\bf Standing notation:}
Throughout this paper we use the following notation:
\begin{itemize}
\item
Given a diffeomorphism $f$ and
a hyperbolic set $\La_f$ of
$f$ there is a $C^1$-neigh\-bor\-hood $\cU$ of $f$ such that
every $g\in \cU$ has a hyperbolic set $\La_g$ called the {\emph{continuation}} of $\La_f$. The set $\La_g$
is close to $\La_f$ and the restrictions of $f$ to $\La_f$ and of $g$ to $\La_g$ are
conjugate. If $P_f$ a hyperbolic periodic point, we denote by $P_g$ the continuation of $P_f$  for $g$ close to $f$.
\item
Given a periodic point $P$ of $f$ we denote
by $\pi(P)$ its period.
\item
The perturbations we consider are always arbitrarily small. Thus
the sentence \emph{there is a $C^r$-perturbation $g$ of $f$}
means
\emph{there is $g$ arbitrarily $C^r$-close to $f$}.
\end{itemize}

\section{Blender-horseshoes}
\label{s.bh}

In this section, we introduce precisely the definition of a 
{\emph{blender-horseshoe,}} a particular case of the blenders in 
\cite{BoDiAnMa}. In fact, blender-horseshoes are the main ingredient
of this paper and the key tool for getting robust homoclinic tangencies.
We beging by reviewing the notion of a blender.

\subsection{Blenders}
\label{ss.blenders} The notion of a $\cut$-blender was introduced
in \cite{BoDiAnMa} as a class of examples, without  a precise and formal
definition.   Blenders were used to get $C^1$-robust
transitivity, \cite{BoDiAnMa}, and robust heterodimensional
cycles, \cite{BDijmj}. The relevance of blenders comes from their
internal geometry and not from their dynamics: a $\cut$-blender is
a (uniformly) hyperbolic transitive set whose stable set robustly
has Hausdorff dimension greater than its stable bundle. In some sense, this property 
resembles and plays a similar role as 
the thick horseshoes introduced by Newhouse, \cite{NeProvidence}.
Following \cite[Definition 6.11]{BDVbook}, we now give a tentative formal definition of a $\cut$-blender:

\begin{defi}[$\cut$-blender]
\label{d.blender} Let $f\colon M\to M$ be a diffeomorphism. A
transitive  hyperbolic set $\Ga$ of $f$ with $\indu(\Ga)=k\ge 2$ is a
{\emph{$\cut$-blender}} if there are a $C^1$-neighborhood $\cU$ of
$f$ and a $C^1$-open set $\cD$ of embeddings of $(k-1)$-dimensional
disks $D$ into $M$ such that, for every diffeomorphism $g\in \cU$,
every disk $D\in \cD$ intersects the local stable manifold
$W^s_{\loc}(\Ga_g)$ of the continuation $\Ga_g$ of $\Ga$ for $g$.
The set $\cD$ is called the {\emph{superposition}} region of the
blender. 
\end{defi}

By definition, the property of a diffeomorphism having a
$\cut$-blender is a  $C^1$-robust property.

We do not know whether $\cut$-blenders yield robust tangencies in the
sense of Definition~\ref{d.rtm}. This leads to the following questions:

\begin{ques}
Let $f\colon M\to M$ be a diffeomorphism  having a $\cut$-blender
$\Ga$ with $k=\indu (\Ga)$.
\begin{itemize}
\item
Does it exist a submanifold $N\subset M$ with  $\dim (N)=k$ such
that that $W^s(\Ga)$ and $N$ have a robust tangency?
\item
Suppose that a submanifold $L$ of dimension $k$  and $W^s(\Ga)$ have
a tangency. Does this tangency yield robust tangencies? More
precisely, does there exist an open set $\cU$ of $\diffM$, $f$ in
the closure of $\cU$,  of diffeomorphisms $g$ with
robust tangencies associated to $\Ga_g$ and ``continuations" of $L$?
\end{itemize}
\label{q.tangency}
\end{ques}

We note that, even for the first
$\cut$-blenders constructed in \cite[Section 1]{BoDiAnMa}, 
these questions remain open. We will give a partial answer to this question in
Theorem~\ref{t.tangency}. For that
we will introduce a special class of $\cut$-blenders, conjugate to
the usual Smale horseshoe,  that we call
{\emph{blender-horseshoes.}}

\subsection{Blender-horseshoes}
\label{ss.blenderhorseshoe}

In this section, we give the precise 
definition of a {\emph{blender-horseshoe}}. This definition involves several concepts
as invariant cone-fields, hyperbolicity, partial hyperbolicity, and Markov partitions,
which we will present separately. Our presentation
follows closely \cite[Section 1]{BoDiAnMa}, thus some details of
our construction are just sketched.

\subsubsection{Cone-fields}\label{sss.conefields}
Consider $\RR^n=\RR^s\oplus\RR\oplus \RR^u$, where $s>0, u>0$, and
$n=s+u+1$. For $\alpha\in(0,1)$, denote by $\cC^{\st}_\alpha$,
$\cC^{\ut}_\alpha $, and $\cC^\uut_\alpha$ the following
cone-fields:
$$
\begin{array}{ll}
\cC^{\st}_\alpha(x)&=\{v=(v^s,v^c,v^u)\in\RR^s\oplus \RR\oplus
\RR^u=T_xM\quad \colon \quad \|v^c+v^u\|\leq \alpha \, \|v^s\|\},\\
 \cC^{\ut}_\alpha(x)&=\{v=(v^s,v^c,v^u)\in\RR^s\oplus \RR\oplus
\RR^u=T_xM\quad \colon \quad
 \|v^s\| \leq \alpha \, \|v^c+v^u\|\},\\
 \cC^{\uut}_\alpha(x)&=\{v=(v^s,v^c,v^u)\in\RR^s\oplus \RR\oplus
\RR^u=T_xM\quad\colon \quad \|v^s+v^c\|\leq \alpha \, \|v^u\|\}.
\end{array}
$$
As $\alpha\in(0,1)$, one gets that $\cC^{\st}_\alpha$ is
transverse to $\cC_\alpha^{\ut}$, that is,
$\cC_\alpha^{\st}(x)\cap \cC_\alpha^{\ut}(x)=0_x\in T_xM$.
Moreover, $\cC_\alpha^{\uut}(x)\subset \cC_\alpha^{\uut}(x)$ for all $x$.

Consider the cube
$$
\CC=[-1,1]^n=[-1,1]^s\times [-1,1]\times [-1,1]^u.
$$
We split the boundary of $\CC$ into three parts:
$$
\begin{array}{ll}
\partial^\st \CC&= \partial \big([-1,1]^s \big) \times [-1,1]\times
[-1,1]^u,\\
\partial^\ct\CC &= [-1,1]^s \times \{-1,1\}\times
[-1,1]^u,
\quad \mbox{and}\\
\partial^\uut \CC&=[-1,1]^s\times [-1,1]\times \partial \big([-1,1]^u\big).
\end{array}
$$
We also consider
$$
\partial^\ut \CC = [-1,1]^s \times \partial \big([-1,1]\times
[-1,1]^u \big)=\partial^c \CC \cup \partial^\uut \CC.
$$
We now consider a local diffeomorphism $f\colon \CC\to \RR^n$ and 
formulate conditions {\bf BH1)}--{\bf BH6)} for the maximal
invariant set $\La$ of $f$ in the cube $\CC$,
$$
\La=\bigcap_{i\in \ZZ}
f^i(\CC),
$$ to be a {\emph{blender-horseshoe,}} see Definition~\ref{d.bd}.

\begin{description}
\item{{\bf BH1)}}\label{ibh1}
The intersection $f(\CC)\cap (\RR^{s}\times \RR \times [-1,1]^u)$ consists of two
connected components, denoted  $f(\cA)$ and $f(\cB)$. Furthermore,
$$f(\cA)\cup f(\cB)\subset (-1,1)^s \times \RR\times [-1,1]^u, \mbox{ and }$$
$$(\cA\cup \cB)\cap \partial^\uut(\CC)=\emptyset.$$
\end{description}

 We denote
$f_{_\cA}\colon \cA\to f(\cA)$ and $f_{_\cB}\colon \cB\to f(\cB)$
the restrictions of $f$ to $\cA$ and $\cB$, respectively. See
Figure~\ref{f.BH1}.

\begin{figure}[htb]

\begin{center}

\psfrag{C}{$\CC$}
\psfrag{Ru}{$\RR^u$}
\psfrag{B}{$\cB$}
\psfrag{A}{$\cA$}
\psfrag{Rc}{$\RR$}
\psfrag{dC}{$\partial^\uut \CC$}

\includegraphics[height=1.6in]{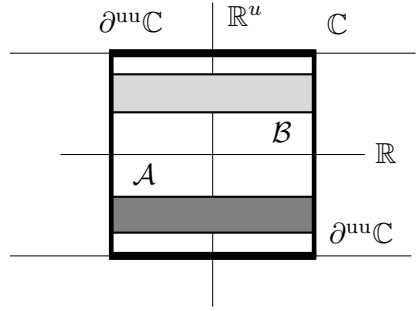}
 \caption{Projection in $\RR\oplus \RR^u$  of a blender horseshoe. Condition {\bf BH1)}.}
\label{f.BH1}

\end{center}

\end{figure}

\begin{description}
\item{{\bf BH2)}} {\bf Cone-fields:}
The cone-field  $\cC_{\alpha}^\st$ is strictly
 $Df^{-1}$-invariant and the cone-fields $\cC^{\ut}_\alpha$ and $\cC^\uut_\alpha$ are strictly
$Df$-invariant. More precisely, there is $0<\alpha'<\alpha$ such
that, for every $x\in  f(\cA)\cup f(\cB)$ one has
$$
Df^{-1}(\cC^\st_\alpha(x))\subset \cC^{\st}_{\alpha'}(f^{-1}x).
$$
In the same way, for every $x\in  \cA\cup \cB$, one has
$$
Df(\cC^\uut_\alpha(x))\subset \cC^{\uut}_{\alpha'}(f(x)) \quad \mbox{and} \quad Df(\cC^\ut_\alpha(x))\subset \cC^{\ut}_{\alpha'}(f(x)).
$$
Moreover, the cone-fields $\cC^\ut_\alpha$ and $\cC^\st_\alpha$ are
uniformly expanding and contracting, respectively.
\end{description}
Note that that property {\bf BH2)} is  open: by increasing slightly
$\alpha'<\alpha$, it holds for every diffeomorphism $g$ in a $C^1$
neighborhood of $f$.

Since
$f(\partial^\uut\CC)$ is disjoint from  $ \RR^{s+1}\times
[-1,1]^u$ and $f(\cA)\cup f(\cB)\subset (-1,1)^s \times \RR\times
[-1,1]^u$, 
from condition
{\bf BH1)} one has that
$$
f(\partial( \CC)) \cap
\partial([-1,1]^{s}\times \RR \times [-1,1]^u) \subset
f(\partial^\st\CC\cup \partial^\ct\CC).
$$ 
Furthermore, this is a $C^1$-robust property.

By {\bf BH2)}, the components of  $f(\partial( \CC)) \cap
\partial([-1,1]^s\times \RR \times [-1,1]^u)$
 are foliated by disks $\Delta$
tangent to $\cC_\alpha^\uut$, i.e., $T_x\Delta \subset \cC^\uut_\alpha(x)$. Hence these disks are transverse to
$\partial([-1,1]^s \times \RR \times [-1,1]^u)$. As a consequence, one gets
the following:
\begin{rema}
Under the ($C^1$-robust) hypothesis {\bf BH2)},  hypothesis { \bf BH1)}
is also a $C^1$-robust property.
\end{rema}

\begin{rema}[Hyperbolicity]
 \label{r.hyperbolicity}
Consider the maximal invariant set $\La$ of $f$ in $\CC$ 
$$
\La=\bigcap_{i\in\ZZ} f^i(\CC).
$$
By {\bf BH1)} and {\bf BH2)} the set $\La$ is  compact and satisfies
$$
\La \subset \mbox{\rm int} (\cA\cup \cB)
\subset \mbox{\rm int} (\CC).
$$
Moreover, the set $\La$ has a dominated splitting $T_\La M=
E\oplus_{_<} F
\oplus_{_<} F$, where $E\subset \cC^\st$, $F\oplus G\subset \cC^\ut$,  and
$G\subset \cC^{\uut}$, and $F$ is one-dimensional.

By {\bf BH2)}, the set $\La$ has a hyperbolic splitting
$E^\st\oplus E^\ut$, where $E^\st=E$ and $E^\ut =F\oplus G$ and
$\dim (E^\st)=s$ and $\dim (E^\ut)=u+1$.
Furthermore, the set $\Lambda$ also has a partially hyperbolic
splitting
$$
T_\La M= E^\st\oplus_{_<} E^{\cut} \oplus_{_<} E^{\uut}
$$
with three non-trivial directions, where $E^{\cut} =F$ and
$E^{\uut}=G$ and
 $\dim (E^\uut)=u$. Note that  $E^\ut=E^\cut \oplus E^\uut$. We say that $E^\uut$ is the {\emph{strong
unstable bundle}} of $\La$.
\end{rema}

\subsubsection{Markov partitions}\label{sss.markov}
Write
$$
\AA=f^{-1}\left(f(\cA)\cap \CC\right)\quad \mbox{and} \quad
\BB=f^{-1}\left(f(\cB)\cap \CC\right).
$$
\begin{description}
\item{{\bf BH3)}}\label{i.ch3} {\bf Associated Markov partition:}

\begin{itemize}
\item
The sets $\AA$ and $\BB$ are both non-empty and connected. That
is, the sets $\AA$ and $\BB$ are the connected components of
$f^{-1}(\CC)\cap \CC$.
\item
The sets $\AA$ and $\BB$ are \emph{horizontal sub-cubes} of $\CC$
and their images $f(\AA)$ and $f(\BB)$ are \emph{vertical
sub-cubes} of $\CC$. More precisely,
$$
f(\AA)\cup f(\BB)\subset  (-1,1)^s \times
[-1,1]\times[-1,1]^u,\mbox{ and}
$$
$$
\AA\cup \BB\subset  [-1,1]^s \times (-1,1)\times(-1,1)^u.
$$
In other words, $f(\AA)\cup f(\BB)$ is disjoint from $\partial^\st
\CC$ and $\AA\cup \BB$ is disjoint from $\partial^\ut\CC$.
\end{itemize}
\end{description}

As a consequence of {\bf BH2)} and {\bf BH3)}, one gets
that $\{\AA,\BB\}$ is a Markov partition generating
$\La$.  Therefore the dynamics of $f$ in $\La$ is
conjugate to the full shift of two symbols. In particular, the hyperbolic
set
$\La$ contains exactly two fixed points of $f$, $P\in \AA$ and $Q\in \BB$. See
Figure~\ref{f.B4}.

\begin{figure}[htb]

\begin{center}

\psfrag{C}{$\CC$}
\psfrag{b}{$\BB$}
\psfrag{B}{$\cB$}
\psfrag{A}{$\cA$}
\psfrag{a}{$\AA$}
\includegraphics[height=1.6in]{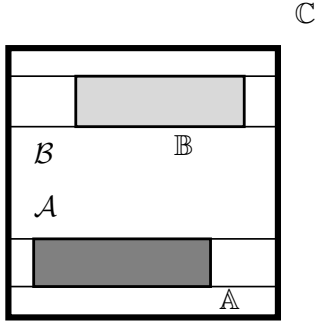}
 \caption{Projection in $\RR \times \RR^u$. {\bf BH3)}
Markov partition of a blender-horseshoe.}
\label{f.B4}

\end{center}

\end{figure}

\subsubsection{$\uut$-disks and their iterates}

\begin{defi}[$\st$- and $\uut$-disks]
A disk  $\De$ of dimension $s$  contained in $\CC$ is an
{\emph{$\st$-disk}} if
\begin{itemize}
\item
it is tangent to $\cC_\alpha^\st$, i.e.,  $T_x \De\subset
\cC_\alpha^\st(x)$ for all $x\in \De$, and
\item
its boundary $\partial \De$ is contained in $\partial^\st (\CC)$.
\end{itemize}
A disk $\Upsilon\subset \RR^s\times\RR\times [-1,1]^u$ of
dimension $u$ is a {\emph{$\uut$-disk}} if
\begin{itemize}

\item it is tangent
to $\cC_\alpha^\uut$, i.e.  $T_x \Upsilon\subset
\cC_\alpha^\uut(x)$ for all $x\in \Upsilon$, and

\item $\partial \Upsilon \subset \RR^s\times\RR\times
\partial([-1,1]^u)$.
\end{itemize}
\end{defi}

Given a point $x\in \La$, there is a unique $f$-invariant manifold
of dimension $u$ tangent at $x$ to the strong unstable bundle
$E^\uut(x)$, the {\emph{strong unstable manifold $W^{\uut}(x)$ of
$x$.}}  For points $x\in \La$, the local invariant manifolds
$W^\st_\loc(x)$, $W^{\ut}_\loc(x)$, and $W^{\uut}_\loc(x)$  are
the connected components of the intersections $W^\st(x)\cap \CC$,
$W^{\ut}(x)\cap \CC$, and $W^{\uut}(x)\cap \CC$ containing $x$,
respectively.

As a consequence of {\bf BH1)--BH3)} one gets

\begin{rema}
For every $x\in\La$, $W^\st_\loc(x)$ is an $\st$-disk and
$W^\uut_\loc(x)$ and is a $\uut$-disk.
\end{rema}

\begin{description}
\item{{\bf BH4)}} {\bf $\uut$-disks through the local stable manifolds of $P$ and $Q$:}
Let $D$ and $D^\prime$ be $\uut$-disks  such that $D\cap
W^\st_\loc(P)\ne \emptyset$ and $D^\prime\cap W^\st_\loc(Q)\ne
\emptyset$. Then
$$
D\cap \partial^\ct(\CC)=D^\prime\cap
\partial^\ct(\CC)=D\cap D^\prime=\emptyset,
$$
%The local stable manifolds $W^\st_\loc(P)$ and $W^\st_\loc(Q)$ are
% $\st$-disks and
% the local strong unstable manifolds $W^\uut_\loc(P)$
%and $W^\uut_\loc(Q)$ are  $\uut$-disks. Moreover,
%$W^\st_\loc(P)\cap W^\uut_\loc(Q)=\emptyset$ and
%$W^\st_\loc(Q)\cap W^\uut_\loc(P)=\emptyset$.
see Figure~\ref{f.B5}.
\end{description}

\begin{figure}[htb]

\begin{center}

\psfrag{D}{$D$}
\psfrag{Dp}{$D^\prime$}
\psfrag{b}{$\BB$}
\psfrag{a}{$\AA$}
\psfrag{c}{$\CC$}
\psfrag{p}{$P$}
\psfrag{q}{$Q$}
\psfrag{wuuP}{$W^\uut_\loc(P)$}
\psfrag{wuuQ}{$W^\uut_\loc(Q)$}
\psfrag{Cuu}{$\cC^{\uut}_\alpha$}
\includegraphics[height=2.0in]{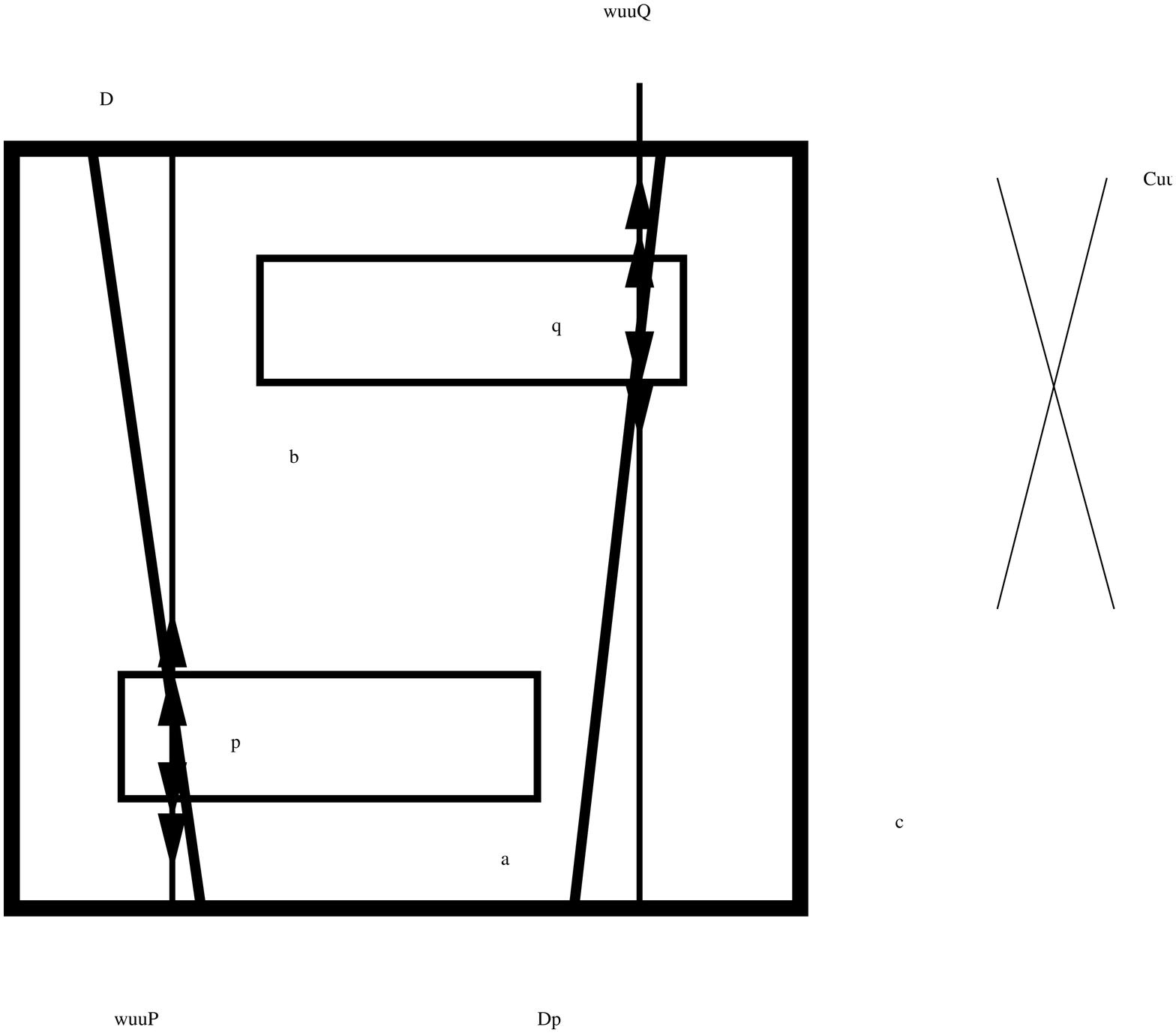}
 \caption{Projection in $\RR\oplus \RR^u$. $\uut$-disks. Condition {\bf BH4)}.}
\label{f.B5}

\end{center}

\end{figure}

Given any  $\st$-disk $\De$, there are two different homotopy
classes of $\uut$-disks contained in $ [-1,1]^s\times
\RR\times[-1,1]^u$ and disjoint from $\De$. We call these classes
{\emph{$\uut$-disks at the right}} and {\emph{at the left}} of
$\De$. We use the following criterion: the $\uut$-disks disjoint
from $W^\st_\loc(P)$ in the homotopy class of $W^\uut_\loc(Q)$ are
at the right of $W^\st_\loc(P)$. The $\uut$-disks disjoint from
$W^\st_\loc(P)$ in the other homotopy class are at the left of
$W^\st_\loc(P)$. We define similarly $\uut$-disks at the left and
at the right of $W^\st_\loc(Q)$, where $\uut$-disks at the left of
$W^\st_\loc(Q)$ are those in the class of $W^{\uut}_\loc (P)$.

According to {\bf BH4)}, $\uut$-disks at the left of
$W^\st_\loc(P)$ are also at the left of $W^\st_\loc(Q)$.
Analogously, $\uut$-disks at the right of $W^\st_\loc(Q)$ are also
at the right of $W^\st_\loc(P)$.

Summarizing, there are five possibilities for a $\uut$-disk $D$
in $ [-1,1]^s\times \RR\times[-1,1]^u$:
\begin{itemize}
\item either $D$ is at the left of $W^\st_\loc(P)$,
\item or $D\cap W^\st_\loc(P)\neq\emptyset$,
\item or $D$ is at the right of $W^\st_\loc(Q)$,
\item or $D\cap W^\st_\loc(Q)\neq\emptyset$,
\item or else $D$ is at the right of $W^{\st}_\loc(P)$
and at the left of $W^{\st}_\loc(Q)$. In this case, we say that the $\uut$-disk $D$ is {\emph{in
between $W^{\st}_\loc(P)$ and $W^{\st}_\loc(Q)$}}.
\end{itemize}

As a consequence of {\bf BH4)} one gets the following.
\begin{rema}[$\uut$-disks in
between $W^{\st}_\loc(P)$ and $W^{\st}_\loc(Q)$] $\,$
\begin{enumerate}
\item
There is a non-empty open
subset $U$ of $\CC$ such that any $\uut$-disk through a point
$x\in U$ is in between $W^{\st}_\loc(P)$ and $W^{\st}_\loc(Q)$. In particular,
there exist $\uut$-disks in between $W^{\st}_\loc(P)$ and
$W^{\st}_\loc(Q)$. 
\item
Every $\uut$-disk $\De\subset[-1,1]^s\times\RR\times[-1,1]^u$ in
between $W^{\st}_\loc(P)$ and $W^{\st}_\loc(Q)$ is contained in
$\CC$  and is disjoint from $\partial^\ct(\CC)$.
\end{enumerate}
\end{rema}

Consider a $\uut$-disk $\De\subset \CC$ and write
$$
f_{_\cA}(\De)=f(\De\cap \cA)
\quad \mbox{and} \quad
f_{_\cB}(\De)=f(\De\cap \cB).
$$
According to {\bf BH1)} and {\bf BH2)} one gets,

\begin{rema}\label{r.uuimages}
For every $\uut$-disk $\De\subset \CC$,
$f_{_\cA}(\De)$ and $f_{_\cB}(\De)$ are
$\uut$-disks in $[-1,1]^s\times \RR\times[-1,1]^u$.
\end{rema}

\begin{description}
\item{{\bf BH5)}} {\bf Positions of images of $\uut$-disks (I):}
Given any  $\uut$-disk $\De\subset \CC$, the following  holds:
\begin{enumerate}
\item
\label{it.1}
if $\De$ is
 at the right of
$W^\st_\loc (P)$ then $f_{_\cA}(\De)$ is a $\uut$-disk at the the
right of $W^\st_\loc (P)$,
\item
\label{it.2}
if $\De$ is
 at the left of
$W^\st_\loc (P)$ then $f_{_\cA}(\De)$  is a  $\uut$-disk at the the
left of $W^\st_\loc (P)$,
\item
\label{it.3}
if $\De$ is
 at the right of
$W^\st_\loc (Q)$ then $f_{_\cB}(\De)$ is a  $\uut$-disk at the the
right of $W^\st_\loc (Q)$,
\item
\label{it.4}
if $\De$ is
 at the left of
$W^\st_\loc (Q)$ then $f_{_\cB}(\De)$ is a  $\uut$-disk at the the
left of $W^\st_\loc (Q)$,
\item
\label{it.5}
if $\De$ is
 at the left of
$W^\st_\loc (P)$ or $\De \cap W^\st_\loc (P) \ne\emptyset$ then
$f_{_\cB}(\De)$ is a  $\uut$-disk at the the left of $W^\st_\loc
(P)$, and
\item
\label{it.6}
if $\De$ is
 at the right of
$W^\st_\loc (Q)$ or $\De \cap W^s_\loc(Q)\ne\emptyset$ then
$f_{_\cA}(\De)$ is a $\uut$-disk at the the right of $W^\st_\loc
(Q)$.
\end{enumerate}
\end{description}
Finally, we state the last condition (which will play a key role) in the definition of
blender-horseshoe:
\begin{description}
\item{{\bf BH6)}} {\bf Positions of images of $\uut$-disks (II):}
Let $\De$ be a  $\uut$-disk in between $W^{\st}_\loc(P)$ and
$W^{\st}_\loc(Q)$. Then either $f_{_\cA}(\De)$ or
$f_{_\cB}(\De)$ is a $\uut$-disk in between
$W^{\st}_\loc(P)$ and $W^{\st}_\loc(Q)$.
\end{description}
Conditions {\bf BH5)}--{\bf BH6)}  are depicted in
Figure~\ref{f.B67}.

\begin{figure}[htb]

\begin{center}

\psfrag{D}{$D$}
\psfrag{E}{$E$}
\psfrag{F}{$F$}
\psfrag{a}{$\AA$}
\psfrag{b}{$\BB$}
\psfrag{wuuP}{$W^\uut_\loc(P)$}
\psfrag{wuuQ}{$W^\uut_\loc(Q)$}
\psfrag{fad}{$f_{_\cA}(D)$}
\psfrag{fae}{$f_{_\cA}(E)$}
\psfrag{fbe}{$f_{_\cB}(E)$}
\psfrag{fbf}{$f_{_\cB}(F)$}

\psfrag{fbd}{$f_{_\cB}(D)$}
\psfrag{faf}{$f_{_\cA}(F)$}

\psfrag{p}{$P$}
\psfrag{q}{$Q$}

\psfrag{text}{$f_{_\cB}(D)$ and $f_{_\cA}(F)$ are not in between $W^\st_\loc(P)$ and $W^\st_\loc(Q)$}

\includegraphics[height=2.1in]{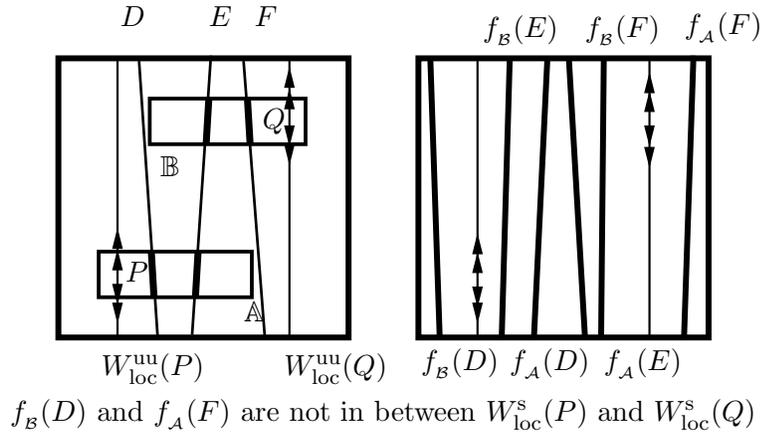}
 \caption{Projection in $\RR\oplus\RR^u$.
Disks in between $W^\st_\loc(P)$ and $W^\st_\loc(Q)$ and their images.}
\label{f.B67}

\end{center}

\end{figure}

%Next definition is a modification of the definition os
%$\cut$-blender in \cite[pages 365-6]{BoDiAnMa}:

\subsubsection{Definition of blender-horseshoe}\label{sss.blenderhorseshoe}
We are now ready for defining  blender-horseshoes\footnote{In some cases, we will
use the terminology $\cut$-blender-horseshoe for
emphasizing that the central one-dimensional direction is
expanding. Using this terminology, a $\cst$-blender-horseshoe is a
blender-horseshoe for $f^{-1}$.}:

\begin{defi}[Blender-horseshoe]
Consider a manifold $M$ of dimension $n\ge 3$ and a diffeomorphism
$f\colon M\to M$. A hyperbolic set $\La$ of  $f$ is a
{\emph{blender-horseshoe}} if (in some coordinate system) there
are a cube $\CC$ and families of cone-fields $\cC^\st$, $\cC^\ut$,
and $\cC^{\uut}$ verifying conditions {\bf BH1)}--{\bf BH6
)}
above.

We say that $\CC$ is the {\emph{reference cube}} of the
blender-horseshoe $\La$ and that the saddles $P$ and $Q$ are the
{\emph{reference
 saddles}}
of $\La$, where {\emph{$P$ is the left saddle}} and {\emph{ $Q$ is
the right saddle.}}
 \label{d.bd}
\end{defi}

%\begin{rema}
%\label{r.complete}
%If $\La$ is a blender horseshoe of a diffeomorphism $f$ then
%$f_{|\La}\colon \La\to \La$ is conjugate to the complete shift of
%two symbols $\sigma\colon \{0,1\}^\ZZ\to \{0,1\}^\ZZ$.
%\end{rema}

%We remark that properties {\bf (BH4-5)} are not satisfied by the
%general blender introduced in \cite{BoDiAnMa}. These properties
%are essential in the proof of Lemma~\ref{l.image}.

Recall that given a hyperbolic set $\La$ of a diffeomorphism $f$
there is a $C^1$-neighborhood $\cU_f$ of $f$ such that every
diffeomorphism $g\in \cU_f$ has a hyperbolic set $\La_g$ which is
close and conjugate to $\La$,  called
the {\emph{continuation}} of $\La$ for $g$. Following \cite[Lemma
1.11]{BoDiAnMa}, one can prove the following:

\begin{lemm}\label{l.continuation}
Let  $\La$ be a blender-horseshoe of a diffeomorphism $f$ with
reference cube $\CC$ and reference saddles $P$ and $Q$. Then there
is a neighborhood $\cU_f$ of $f$ in $\diffM$ such that for all
$g\in\cU_g$ the continuation $\La_g$ of $\La$ for $g$ is a
blender-horseshoe with reference cube $\CC$ and reference saddles
$P_g$ and $Q_g$.
\end{lemm}

Arguing  as in \cite{BoDiAnMa} (see also \cite{BDV} for a simple toy
model and the example in
Section~\ref{ss.prototype}) one gets the following:

\begin{rema}
Every $\uut$-disk in between in between $W^{\st}_\loc(P)$ and
$W^{\st}_\loc(Q)$ intersects $W^\st_\loc(\La)$.
%Hence the blender-horseshoe $\La$ is a $\cut$-blender in the sense of
%Definition~\ref{d.blender}):
Therefore the blender-horseshoe $\La$ is a $\cut$-blender in the
sense of Definition~\ref{d.blender}, where  the $\uut$-disks in
between $W^{\st}_\loc(P)$ and $W^{\st}_\loc(Q)$ define its
superposition region.
\end{rema}

\section{Robust tangencies}
\label{s.robusttangencies} In this section, we introduce a class
of sub-manifolds called folding manifolds relative to a
blender-horseshoe $\La$. The main technical step is
Proposition~\ref{p.tangency}, which claims that folding manifolds are
tangent to the local stable manifold  $W^\st_\loc(\La)$ of the blender-horseshoe
$\La$. Moreover, these tangencies are $C^1$-robust, see
Theorem~\ref{t.tangency}. Finally,  using blender-horseshoes, in Theorem~\ref{th.zeghib} we give sufficient conditions 
for the generation of
robust tangencies by 
a homoclinic
tangency.

\subsection{Folding manifolds and tangencies associated to
blender-horseshoes}\label{s.tangencies}

Let $\La$ be a blender-horseshoe
with reference cube $\CC$
 as in  Section~\ref{ss.blenderhorseshoe}. Recall that the dimension of the unstable bundle of $\La$, 
 $E^\ut=E^\cut \oplus E^\uut$, is
$(u+1)$. We say that \emph{the $\ut$-index of the
blender-horseshoe $\La$ is $(u+1)$} (i.e., the $\ut$-index of
$\La$ as a hyperbolic set).
Define the {\emph{local stable manifold of $\La$}} by
$$
W^\st_{\loc}(\La)=\bigcup_{x\in \La} W^\st_{\loc}(x),
$$
where
$W^\st_{\loc}(x)$ is the connected component of $W^\st(x)\cap \CC$
containing $x$. .

\begin{rema}\label{r.stableset}
The local stable manifold $W^\st_{\loc}(\La)$ of the
blender-horseshoe $\La$ is the set of points $x\in \CC$ whose
forward orbit remains in the reference cube $\CC$.
\end{rema}

A new ingredient of this section is the notion of folding manifold defined as follows:

\begin{defi}[Folding manifold]
 \label{d.folding}
 Consider a blender-horseshoe $\Lambda$ of $\ut$-index $(u+1)$ with reference cube $\CC$
 and reference saddles $P$ and $Q$.
A sub-manifold $\cS\subset \CC$ of dimension $(u+1)$ is a
{\emph{folding manifold of $\La$ (relative to the saddle $P$)}} if there is
a  family $(\cS_t)_{t\in [0,1]}$ of $\uut$-disks
depending continuously on $t$ such that:
\begin{itemize}
\item
$\cS=\bigcup_{t\in [0,1]} \cS_t,$
 \item $\cS_0$ and $\cS_1$ intersects $W^\st_\loc (P)$, and
 \item for every $t\in (0,1)$, the $\uut$-disk $\cS_t$ is in
 between $W^\st_\loc(P)$ and $W^\st_\loc(Q)$.
\end{itemize}
We similarly define a {\emph{folding manifold (relative to $Q$)}}.
A {\emph{folding manifold}} of the blender-horseshoe $\La$  is a
folding manifold relative either to $P$ or to $Q$.
\end{defi}

\begin{rema}\label{r.folding}
 Let $\cS$ be a folding manifold of the blender-horseshoe  $\La$. Then there are a
point $x\in \cS$ and a non-zero vector $v\in T_x \cS$ such that
$v\in \cC^\st(x)$.
\end{rema}

A key property of 
folding manifolds of blender-horseshoes is the following:

\begin{prop}
\label{p.tangency} Let $\cS$ be  a folding manifold  of a
blender-horseshoe $\La$. Then $\cS$ and $W^\st_\loc(\La)$ are
tangent at some point $z$.
\end{prop}

To prove this proposition we need the following lemma:

\begin{lemm}
\label{l.image} Consider a diffeomorphism $f$ having a
blender-horseshoe $\La$ as above. The image by $f$ of a folding
manifold $\cS$ of $\La$ contains a folding manifold of $\La$.
\end{lemm}

\begin{demo}
Let us assume, for instance, that the folding manifold $\cS$ is relative to $P$.
We will prove that either $f_{_\cA}(\cS)$ is a folding manifold
relative to $P$ or $f_{_\cB}(\cS)$ contains a folding manifold
relative to $P$.

If $f_{_\cA}(\cS)$ is a folding manifold relative to $P$ we are
done. So we can assume that $f_{_\cA}(\cS)$ is not a folding
manifold. As the $\uut$-disks $\cS_0$ and $\cS_1$ meet
$W^\st_\loc(P)$, by Remark~\ref{r.uuimages}, their images
$f_{_\cA}(\cS_0)$ and $f_{_\cA}(\cS_1)$ are $\uut$-disks
intersecting $W^\st_\loc(P)$. Furthermore, by item~\ref{it.1} in
{{\bf BH6)}}, the $\uut$-disk $f_{_\cA}(\cS_t)$ is at the right of
$W^\st_\loc(P)$, for every $t\in(0,1)$.

 Since we are assuming that
$f_{_\cA}(\cS)$ is not a folding manifold relative to $P$, by
definition of a folding manifold relative to $P$,
 there is some $t_0\in (0,1)$ such that
$f_{_\cA}(\cS_{t_0})$ is not at the left of $W^\st_\loc(Q)$ (i.e., it is either at 
the left of $W^\st_\loc(Q)$ or it meets $W^\st_\loc(Q)$).
Thus, by
continuity of the disks $f_{_\cA}(\cS_t)$, there is $t_1\in (0,t_0)$ such that $f_{_\cA}(\cS_{t_1})
\cap W^\st_\loc (Q)\ne \emptyset$. As $\cS_{t_1}$ is in
 between $W^\st_\loc(P)$ and $W^\st_\loc(Q)$, by {\bf BH6)}, the image $f_{_\cB}(\cS_{t_1})$ is in between
 $W^\st_\loc (P)$ and $W^\st_\loc(Q)$.

By the definition of folding manifold relative to $P$, every
$\cS_t$ is at the left of $W^\st_\loc(Q)$. Therefore, by
item~\ref{it.4} in {\bf BH5)}, $f_{_\cB}(\cS_t)$ is a $\uut$-disk
at the left of $W^\st_\loc (Q)$, for every $t\in (0,1)$. Moreover,
by item~\ref{it.5} in {\bf BH5)},  the images $f_{_\cB}(\cS_0)$
and $f_{_\cB}(\cS_1)$ are $\uut$-disks at the left of
$W^\st_{\loc}(P)$.

Now by continuity of the disks $f_{_\cB}(\cS_t)$ and since
$f_{_\cB}(\cS_{t_1})$ is at the right of $W^\st_\loc (P)$), there
are
 parameters $t_2$ and $t_3$, with $t_2<t_1<t_3$, such that
$f_{_\cB}(\cS_{t_2})$ and $f_{_\cB}(\cS_{t_3})$ are $\uut$-disks
intersecting $W^\st_{\loc}(P)$ and  $f_{_\cB}(\cS_{t})$ is a
$\uut$-disk at the right of $W^\st_\loc (P)$, for all $t\in
(t_2,t_3)$. Since, by item~\ref{it.4} in {\bf BH5)}, these disks
are at the left of $W^\st_\loc (Q)$, they are in between
$W^\st_\loc (P)$ and $W^\st_\loc (Q)$. This implies that
$$
\bigcup_{t\in [t_1,t_2]} f_{_\cB}(\cS_t) \subset f_{_\cB}(\cS)
$$
is a folding manifold relative to $P$, ending
the proof of the lemma.
\end{demo}

We are now ready to conclude the proof of Proposition~\ref{p.tangency}.

\begin{demo}[Proof of Proposition~\ref{p.tangency}]
Write $\cS_0=\cS$. By Lemma~\ref{l.image}, there is a folding
manifold $\cS_1$ contained in $f(\cS_0)$. Using
Lemma~\ref{l.image} and arguing inductively, we define a sequence
of folding manifolds $(\cS_i)_i$ of the blender-horseshoe $\La$ such that, for every
$i\ge 0$, $\cS_{i+1}$ is contained in $f(S_i)$. Let 
$$
\tilde
\cS_i=f^{-i}(S_i).
$$ 
In this way we get a nested sequence $(\tilde
\cS_i)_i$, $\tilde \cS_{i+1}\subset \tilde \cS_i\subset \cS$, of
connected and compact sets. Thus, by construction, the intersection set
$$
\cS_\infty=\bigcap_{i=0}^\infty \tilde\cS_i \ne \emptyset
$$ 
is  connected and compact. Moreover, $\cS_\infty\subset \cS_0$.

 By construction, the whole forward
orbit of the set $\cS_\infty$ is contained in the reference cube $\CC$ of
the blender. By Remark~\ref{r.stableset}, the set $\cS_\infty$ is
contained in
 $W^\st_\loc(\La)$. Note that, as $\La$ is totally
disconnected (a Cantor set) and $\cS_\infty$ is connected, hence there is
some $z\in \La$ such that 
$$
\cS_\infty\subset W^\st_\loc (z).
$$
Since $\cS_i$ is a folding manifold for every $i$, 
Remark~\ref{r.folding} implies that there are a point $x_i\in \cS_i$ and a
non-zero vector $v_i\in T_{x_i}\cS_i$ such that
$v_i\in\cC^\st_\alpha (x_i)$. 

Consider the point $\tilde
x_i=f^{-i}(x_i)\in \tilde\cS_i\subset \cS$ and an unitary vector
$\tilde v_i$ parallel to $D_{x_i} f^{-i} (v_i)$. Note that $\tilde
v_i \in T_{\tilde x_i} \tilde \cS_i$, thus $\tilde v_i \in
T_{\tilde x_i} \cS$. By the $(Df^{-1})$-invariance of the cone-field
$\cC^\st_\alpha$, condition {\bf BH2)}, we have that $\tilde v_i
\in \cC^\st_\alpha(\tilde x_i)$.
 We can assume (taking a subsequence if necessary) that
 $$
 \tilde
 x_i\to  x_\infty\in \cS_\infty\subset \cS\cap W^\st_\loc(\La)
 \qquad \mbox{and} \qquad \tilde v_i\to v_\infty.
 $$
Hence $x_\infty\in \cS_\infty$ and $x_\infty\in W^\st_\loc(z)$.
Our construction also implies that $v_\infty \in T_{x_\infty}
\cS$. Finally, also by construction, the vector $v_\infty$ belongs
to the intersection
$$\bigcap_{i\ge 0}
Df^{-i} (\cC^\st_\alpha
(f^i(x_\infty)))=T_{x_\infty}W^\st_\loc(z).$$
 This completes the
proof of the proposition.
\end{demo}

\subsection{Robust tangencies}
\label{ss.newrobust}

We now return to the problem of robust tangencies in Question~\ref{q.tangency}. 
We need the following definition.

\begin{defi}[Folded  submanifolds with respect to a blender-horseshoe]
Let $f\colon M\to M$ be a diffeomorphism having a
blender-horseshoe $\Lambda$ with reference saddles $P$ and $Q$ and $N\subset M$ be a submanifold  of
dimension $\indu (\La)$.

We say that $N$ is {\emph{folded
with respect to $\La$}} if the interior of $N$ contains a folding
manifold $\cS=(\cS_t)_{t\in [0,1]}$ relative to some reference
saddle $A\in\{P,Q\}$ of the blender, here $(\cS_t)_{t\in [0,1]}$ is the family of $\uut$-disks in
Definition~\ref{d.folding}, such that:
\begin{itemize}
\item $\cS_0\cap W^\st_\loc(A)$ and $\cS_1\cap W^\st_\loc(A)$ are
transverse intersection points of $N$ with $W^\st_\loc(A)$;
\item
There is $0<\alpha'<\alpha$ such that  the $\uut$-disks $\cS_t$,
$t\in[0,1]$, are tangent to the cone-field $\cC^\uut_{\alpha'}$.
\end{itemize}
To emphasize the reference saddle $A$ of the blender
we consider, we say
that the submanifold $N$ is {\emph{folded with respect to
$(\La,A)$.}}
\end{defi}

\begin{rema}\label{r.folded} A submanifold to be folded with respect to a blender-horseshoe is
a $C^1$-open property.
\end{rema}

As a direct consequence of Proposition~\ref{p.tangency} and
Remark~\ref{r.folded} one gets:

\begin{theo}
\label{t.tangency} Let $N\subset M$ be  a  folded submanifold 
with respect to a blender-horseshoe $\La$. Then $N$ and
$W^\st_\loc(\La)$ have a $C^r$-robust tangency.
\end{theo}

\subsection{Robust homoclinic tangencies}\label{ss.robust}

In this section we prove that homoclinic tangencies associated to 
blender-horseshoes yield $C^r$-robust homoclinic tangencies.

\begin{theor}\label{th.zeghib}
Consider a transitive hyperbolic set
$\Si$ of a  $C^r$-diffeomorphism $f$ containing a $\cut$-blender-horseshoe and
a saddle with a homoclinic tangency. Then there is a diffeomorphism $g$
arbitrarily $C^r$-close to $f$ such that the continuation $\Si_g$
of $\,\Si$ has a $C^r$-robust homoclinic tangency.
\end{theor}

We need the following lemma.

 \begin{lemm}\label{l.tangency} Consider a $C^r$-diffeomorphism $f$ with  a $\cut$-blender-horseshoe
 $\La$.
 Assume that there is a saddle $R$ with $\indu (R)=\indu (\La)$ and such
 that $W^\ut(R)$ has a tangency with $W^\st(A)$, where $A$ is a reference saddle of the blender $\La$.
 Then there is a diffeomorphism $g$ arbitrarily $C^r$-close to $f$ such that
 $W^\ut(R_g)$ is a folded manifold with respect to the continuation
 $\La_g$ of the blender-horseshoe $\La$.
 \end{lemm}
 By Theorem~\ref{t.tangency} one gets:

 \begin{coro}\label{c.tangency}
In Lemma~\ref{l.tangency}, the  stable manifold $W^\st_\loc(\La_g)$  of the blender-horseshoe  and
 $W^\ut(R_g)$ have a $C^r$-robust tangency.
 \end{coro}

 \begin{demo}[Proof of Lemma~\ref{l.tangency}]
We suppose that $A=P$ is the left reference saddle of the blender.
The proof  involves a string of $C^r$-perturbations of
the diffeomorphism $f$. For simplicity, we also denote these perturbations by $f$.

We begin by noting that \emph{the center stable bundle $E^{\cst}$}
is well defined for every point $x$ in the local stable manifold
$W^\st_\loc(\La)$. Recall that $W^\st_\loc(\La)$ is the set of points whose forward orbit
remains in the reference cube $\CC$ of the blender-horseshoe, Remark~\ref{r.stableset}.
Given a point $x\in W^\st_\loc(\La)$, the subspace $E^{\cst}(x)$
is the set of vectors $v\in T_xM$ such that $Df^n(v)\notin
(\cC^\uut(f^n(x))\setminus\{0\})$, for every $n\geq 0$. The space
$E^\cst(x)$ has dimension $\inds(\La)+1$ and depends continuously
on the point $x\in W^\st_\loc(\La)$ and on the diffeomorphism $f$.

First, 
after considering forward iterations, we can assume that the
tangency intersection point $B$ between $W^\ut(R)$ and $W^\st(P)$  is in
$W^\st_\loc(P)$. Therefore the whole forward orbit of  $B$ is the
reference cube $\CC$. Recall that $P$ is the left reference saddle of
the blender-horseshoe, thus $P$ is in the ``rectangle'' $\AA$ of the Markov
partition. Thus,  for any $n\geq 0$,  $f_{_\cA}^n(B)$ is defined and
belongs to $\CC$. Hence $E^\cst(f^i(B))$ is well defined for all
$i\geq 0$ (recall the comment above). Note that
$$
\dim(T_B W^\ut (R))+\dim
(E^\cst(B))=\dim(T_B W^\ut (R))+\dim (T_B W^\st(P))+1=
\dim (M)+1.
$$
Thus
after a  perturbation,
we can assume that    $T_B W^\ut (R)$  is transverse
to $E^\cst(B)$.  Hence there is a subspace  $\VV\subset T_B W^\ut
(R)$ of dimension $u$ with $\VV\oplus E^\cst(B)= T_B M$, where $(u+1)=\indu (\La)$.

Consider any $\alpha'\in(0,\alpha)$ ($\alpha$ is the constant in the definition of the cone-fields of
the blender). Then, for every $n>0$ large enough,
one has that $Df^n(\VV)$ is contained in the cone $\cC^\uut_{\alpha^\prime}(f^n(B))$. For
simplicity, let us assume that $n=0$, that is $\VV\subset
\cC^\uut_{\alpha^\prime}(B)$.
This implies that (up to increase slightly the constant $\alpha'<\alpha$) there is a
small submanifold $\tilde \cS\subset W^\ut (R)$ such that
the point  $B$  is in  the interior of $\tilde\cS$, and 
$\tilde \cS$ is foliated by disks $(\tilde \cS_t)_{t\in [-1,1]}$ of
dimension $u$ tangent to the  cone-field
$\cC^\uut_{\alpha^\prime}$.

Using the expansion by $Df$  in the cone-field $\cC^\uut_{\alpha^\prime}$ and considering forward
iterations of $\tilde \cS$ by $f$,
we get $k\ge 0$ and a submanifold  $\cS\subset f^k(\tilde \cS)$ such that:
\begin{itemize}
 \item $\cS$ contains
$f^k(B)$ in its interior and  is tangent to $W^\st (P)$ at $f^k(B)$, and
\item
$\cS$
is foliated by $\uut$-disks $\cS_t\subset f^k(\tilde \cS_t)$,
$t\in[-1,1]$, where the disks $\cS_t$  are tangent to $\cC^\uut_{\alpha'}$.
\end{itemize}
Again for simplicity, we assume that $k=0$ and that $B\in\cS_0$.

After a new perturbation, we can assume that the
contact between $\cS$ and $W^\st_\loc(P)$  at the point $B$
is quadratic. In particular,
there is  small $\epsilon > 0$, such that either all the
$\uut$-disks $\cS_t$, $t\ne 0$ and $t\in [-\epsilon,\epsilon]$,  are  at
the left of $W^\st_\loc(P)$ (case (a)), or all
the
$\uut$-disks $\cS_t$, $t\ne 0$ and $t\in [-\epsilon,\epsilon]$,  are
at the right of $W^\st_\loc(P)$ (case (b)). So after discarding some disks and
reparametrizing the family $\cS_t$,
 we can assume that $\epsilon=1$.

\begin{figure}[htb]

\begin{center}

\psfrag{p}{$P$}
\psfrag{s}{$\cS$}
\psfrag{text}{perturbation}
\psfrag{proj}{projection in $\RR^s\oplus \RR^u$}

\includegraphics[height=2.0in]{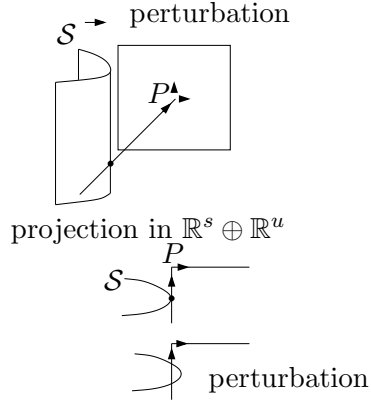}
 \caption{Folded manifold: case (a).}
\label{f.casoa}

\end{center}

\end{figure}

\medskip

{\bf Case (a):}
{\emph{for $t\ne 0$, every $\cS_t$  is at the left of $W^\st_\loc(P)$.}}
After a new perturbation, we can assume  that $\cS$ is a folding manifold
relative to $P$, see Figure~\ref{f.casoa}. Since, by construction,  $\cS$ is contained in $W^\ut(R)$ this concludes
the proof in the first case.

\medskip

{\bf Case (b):} 
{\emph{for $t\ne 0$ every $\cS_t$  is at the right of $W^\st_\loc(P)$.}}
By considering positive iterations of $\cS$ by $f_{_\cA}$,
one gets a large $i>0$ such that $f_{_\cA}^i(\cS)$ meets transversely
$W^\st_\loc(Q)$ at some points in $f_{_\cA}^i(\cS_{t_1})$ and
$f_{_\cA}^i(\cS_{t_2})$, where $t_1<0<t_2$. Once again, let us assume that $i=0$.

\begin{figure}[htb]

\begin{center}

\psfrag{p}{$P$}
\psfrag{q}{$Q$}
\psfrag{s}{$\cS$}
\psfrag{q}{$Q$}
\psfrag{fi}{$f^i$}
\psfrag{fs}{$f_{_\cA}^i(\cS)$}
\psfrag{text}{perturbation}
\psfrag{proj}{projection in $\RR^s\oplus \RR^u$}
\includegraphics[height=2.4in]{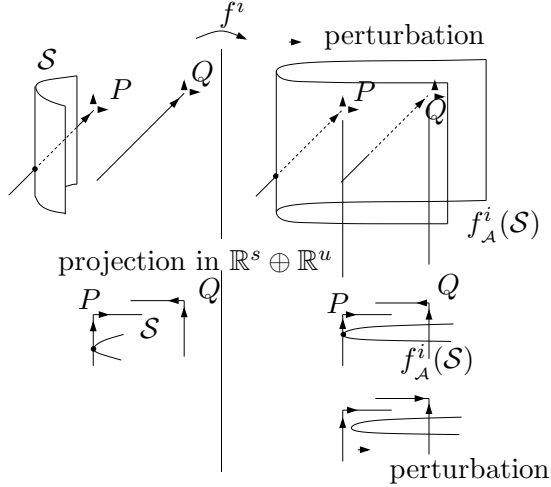}
 \caption{Folded manifold: case (b).}
\label{f.casob}

\end{center}

\end{figure}

We can choose $t_1$ and $t_2$ such that
 the disks $\cS_t$ are at the left of $W^s_\loc(Q)$ for every  $t\in (t_1,t_2)$. Recall that, by hypothesis, for all $t\in[t_1,t_2]\setminus\{0\}$
the  disks
$\cS_t$ are at the right of
$W^\st_\loc(P)$. Therefore,  we can perform a final perturbation so that
$\hat \cS=\bigcup_{t\in[t_1,t_2]}\cS_t$ is a folding manifold
relative to $Q$, see Figure~\ref{f.casob}. Since $\hat \cS$ is
contained in $W^\ut(R)$ this ends the proof of the lemma.
\end{demo}

We are now ready to prove Theorem~\ref{th.zeghib}.

\begin{demo}[Proof of Theorem~\ref{th.zeghib}] It is enough to
observe that after a $C^r$-perturbation, one can assume that
the homoclinic tangency of $\Sigma$ occurs
between the unstable manifold of a periodic point $R\in\Si$ and
the left reference saddle of the blender-horseshoe.
\end{demo}

\section{Generation of blender-horseshoes}
\label{s.examples}

In this section we see  how blender-horseshoes arise naturally in our non-hyperbolic setting.
First, in Section~\ref{ss.prototype}, we review constructions in  \cite{BDV} 
providing simple examples of  blender-horseshoes. 
In
Section~\ref{sss.affine}, using these constructions,
we see that partially hyperbolic saddles (saddle-node and flip points) with strong homoclinic
intersections (intersections between the strong stable and strong unstable manifolds of a non-hyperbolic saddle) yield 
blender-horseshoes.
Finally, following
\cite{BDijmj}, in Sections~\ref{sss.coindexblender} and \ref{sss.profofweak} 
we prove that co-index one heterodimensional cycles
generate blender-horseshoes. 
We will see in Section~\ref{ss.generation}
that co-index one cycles occur naturally in the non-hyperbolic setting.

\subsection{Prototypical blender-horseshoes}
\label{ss.prototype}
In this section,
we
consider a
local diffeomorphism $f$ having an affine horseshoe $\La$ with a
dominated splitting
with three non-trivial bundles,  $E^{\st}\oplus_{_<} E^\cut\oplus_{_<} E^\uut$, where $E^\cut$ is one dimensional
and $E^\ut=E^\cut \oplus E^\uut$ is the unstable bundle of $\La$.
We suppose that
 $\La$ is contained in a hyperplane  $\Pi$
tangent to $E^\st\oplus E^{\uut}$ and that the expansion along the
direction $E^\cut$ is close to one. Under these assumptions we prove that there are
perturbations $g$ of $f$ such that the continuations $\La_g$
of $\La$ for $g$ are blender-horseshoes, see Proposition~\ref{p.blenderhorseshoe}. 
We now go to the details of this construction, we borrow from 
\cite{BDV}. 

\medskip

Let  $\DD=[-1,1]^n$ and $n=s+u,\, s,u\ge 1$. Consider a
diffeomorphism $F\colon \RR^n\to \RR^n$  having a horseshoe
$\Sigma=\cap_{k\in \ZZ} F^k(\DD)$ such that:
\begin{itemize}
\item
 $F^{-1}(\DD) \cap
\DD$ consists of two connected components $\DD_1=[-1,1]^s
\times \UU_1$ and $\DD_2=[-1,1]^s\times \UU_2$, where $\UU_1$
and $\UU_2$ are disjoint topological compact disks of dimension
$u$.
\item
The map $F$ is affine on each rectangle $\DD_i$: there are
linear maps $S_i \colon \RR^s\to \RR^s$ and $U_i\colon \RR^u\to
\RR^u$, $i=1,2$, such that
$$
DF|_{\DD_i}= \left(
\begin{matrix}
S_i & 0\\ 0 & U_i
\end{matrix}
\right),\qquad ||S_i||, ||U_i^{-1}||<1/2, \quad i=1,2,
$$
where $||A||$ is the norm of the linear map $A$.
\end{itemize}
 We suppose that,
in the usual coordinates $(x^s,x^u)$ in $\RR^n=\RR^s\times \RR^u$,
 the fixed saddles of the horseshoe $\Si$ are $p=(0^s,0^u)\in 
\DD_1$ and $q=(a^s,a^u)\in \DD_2$.

\begin{figure}[htb]

\begin{center}

\psfrag{Ru}{$\RR^u$}
\psfrag{Rs}{$\RR^s$}
\psfrag{R}{$\RR$}
\psfrag{D1}{$\DD_1$}
\psfrag{D2}{$\DD_2$}
\psfrag{I}{$[-1,1]^n$}

\includegraphics[height=2.0in]{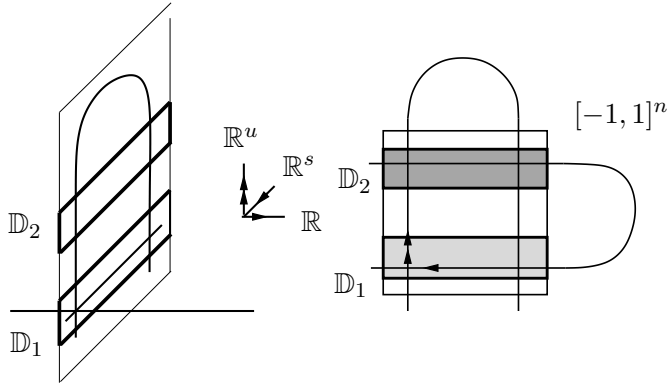}
 \caption{An affine horseshoe. The map $f_{\la,0}$.}
\label{f.horseshoe}

\end{center}

\end{figure}

Consider $\lambda\in (1,2)$ and the family of local
diffeomorphisms $(f_{\la,\mu})_{\mu\in [-\epsilon,\epsilon]}$ of
$\RR^{n+1}$ given by:
$$
f_{\la,\mu}(x^s,x^u,x)=\left\{
\begin{array}{lll}
&(F(x^s,x^u),\la\,x), &\quad  \mbox{if $x^u\in \UU_1$},\\
&(F(x^s,x^u),\la\,x-\mu),&\quad \mbox{if $x^u\in \UU_2$}.
\end{array}
\right.
$$
By definition, for small $\mu$, the diffeomorphism $f_{\la,\mu}$
has two fixed saddles: $P=(0^s,0^u,0)$ (independent of $\la$ and
$\mu$) and $Q_{\la,\mu}=(a^s,a^u,\mu/(\la-1))$.

Let $\La_{\la,0}$ be the maximal invariant set of $f_{\la,0}$ in 
$(\DD_1\cup \DD_2)\times [-1,1]$.
Note that
$\La_{\la,0}=\Sigma\times \{0\}$ is a hyperbolic set of
$f_{\la,0}$. We say that $\La_{\la,0}$ is an {\emph{affine
horseshoe}} of $f_{\la,0}$ with central expansion $\lambda$.
Observe that
that the hyperplane $\RR^n\times \{0\}$ is not normally hyperbolic for $f_{\la,0}$.

 We denote by $\La_{\la,\mu}$ 
the maximal invariant set of $f_{\la,\mu}$ in 
$(\DD_1\cup \DD_2)\times [-1,1]$. For small $\mu$,
the set $\La_{\la,\mu}$  is 
the continuation of
$\La_{\la,0}$. More precisely,
fixed small $\delta>0$ and the cube $\CC_\de=[-1,1]^s\times
[-1,1]^u \times [-\de,\de]$, for $|\mu|<(\la-1)\,\de$, the
set $\La_{\la,\mu}$ is the maximal invariant set of $f_{\la,\mu}$
in $\CC_\de$.
Clearly, $P,Q_{\la,\mu}\in
\La_{\la,\mu}$.

%Moreover,
%$\La_t$ is the homoclinic class of $P$ (or $Q_t$).

\begin{prop}
 \label{p.blenderhorseshoe}
For  every $\lambda>1$ close to $1$ and $\mu>0$, the set
$\La_{\la,\mu}$ is a blender-horseshoe with reference cube
$\CC_\de$ and reference saddles $P$ and $Q_{\la,\mu}$ ($P$ is the
left saddle and $Q_{\la,\mu}$ the right one).
\end{prop}

In this section, for notational convenience, we
write the central coordinates in the third position.

\begin{demo}We fix $\lambda>1$ and $\mu>0$ and we simply write
$\La$, $f$, $P$, and  $Q$,  omitting the dependence on the
parameters. The hyperbolicity of $\La$ follows from the
hyperbolicity of $F$ and from the normal expansion by $\lambda>1$.
Consider the constant bundles
$$
E^\st= \left(\RR^{s}\times\{(0^u,0)\}\right), \quad E^\ct=
\left(\{0^s,0^u\}\times \RR\right), \quad
E^\uut=\left(\{0^s\}\times \RR^u\times\{0\}\right).
$$
Since $\la$ is less than $2$, then
$$
T_x\RR^{n+1}=E^\st\oplus_{_<} E^\ct\oplus_{_<} E^\uut, \quad x\in
\La,
$$
is a dominated splitting over $\La$. Furthermore, as the bundles above
are constant,  the cone-fields $\cC^\cut_\alpha$ and
$\cC^\uut_\alpha$ are $Df$-invariant and $\cC^\st_\alpha$ is
$(Df^{-1})$-invariant, for every $\alpha\in(0,1)$. Finally, for small $\alpha$,
$\cC^\st_\alpha$ is uniformly contracting and $\cC^\ut_\alpha$
is uniformly expanding.
This gives
condition {\bf BH2)}.

To get conditions {\bf BH1)} and {\bf BH3)} just let
$$\cA=[-1,1]^s \times
U^{-1}_1([-1,1]^u)\times [-\de,\de],\qquad
\cB=[-1,1]^s \times
U^{-1}_2([-1,1]^u)\times [-\de,\de],$$
$$
\AA=[-1,1]^s \times
U^{-1}_1([-1,1]^u)\times
[-\lambda^{-1}\, \de,\lambda^{-1}\, \de],\qquad\BB=[-1,1]^s \times
U^{-1}_2([-1,1]^u)\times
\left[\frac{-\de+\mu}\la,\frac{\de+\mu}\la\right],
$$ and observe that
$$f(\AA)=S_1([-1,1]^s) \times
[-1,1]^u\times [-\de,\de]\quad \mbox{ and } \quad f(\BB)= S_2([-1,1]^s)
\times [-1,1]^u\times [-\de,\de].$$

\begin{figure}[htb]

\begin{center}

\psfrag{Cd}{$\CC_\de$}
\psfrag{cB}{$\cB$}
\psfrag{cA}{$\cA$}
\psfrag{BB}{$\BB$}
\psfrag{AA}{$\AA$}
\psfrag{flm}{$f_{\la,\mu}$}

\includegraphics[height=1.8in]{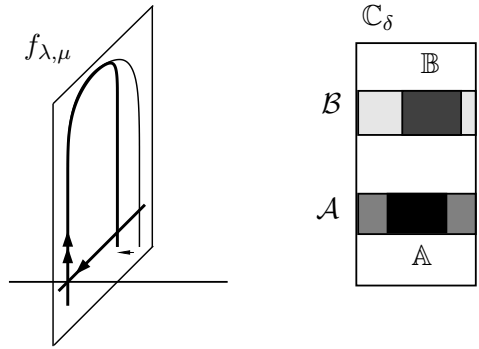}
 \caption{Prototypical blender-horseshoe.}
\label{f.prototypal}

\end{center}

\end{figure}

Observe that the local invariant manifolds of $P$ and $Q$ are:
$$\begin{array}{lll}
&W^{\st}_{\loc}(P,f) =[-1,1]^s\times \{(0^u,0)\},\qquad
&W^{\st}_{\loc}(Q,f) =[-1,1]^s\times \{(a^u, \frac{\mu}{\la-1})\},\\
\\ &W^{\uut}_{\loc}(P,f)= \{0^s\} \times [-1,1]^u \times \{0\},\qquad
&W^{\uut}_{\loc}(Q, f)= \{a^s\} \times [-1,1]^u \times
\{\frac{\mu}{\la-1}\}.
\end{array}
$$
It is immediate to check that  vertical disks of the form
$\{x^s\}\times [-1,1]^u\times\{x^c\}$ satisfy condition {\bf
BH4)}. To get {\bf BH4)} for $\uut$-disks it is enough to take
$\alpha\in(0,1)$ small enough in the definition of the
cone-fields.

Condition {\bf BH5)} follows from the fact that $f_{_\cA}$ and
$f_{_\cB}$ are affine maps preserving the dominated splitting and
whose center eigenvalue $\lambda$ is positive.

\begin{figure}[htb]

\begin{center}

\psfrag{wuuP}{$W^\uut_\loc(P,f)$}
\psfrag{wuuQ}{$W^\uut_\loc(Q,f)$}
\psfrag{de}{$\De$}
\psfrag{I1}{$I_1$}
\psfrag{I2}{$I_2$}

\includegraphics[height=1.8in]{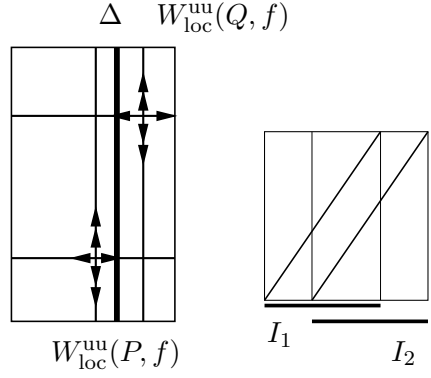}
 \caption{$\uut$-disks and one-dimensional reduction.}
\label{f.onedim}

\end{center}

\end{figure}

It remains to check condition {\bf BH6)}.  We first consider
vertical disks $\De$ parallel to $E^\uut$
$$
\De=\{x^s\} \times [-1,1]^u \times \{x^c\}
$$
which are in
between $W^{\st}_{\loc}(P,f)$ and $W^{\st}_{\loc}(Q,f)$. This
means that $x^c\in (0, \frac{\mu}{\la-1})$. We consider two cases:

\noindent{\bf Case 1:} $x^c\in I_1= (0,\frac{\mu}{\la\, (\la-1)})$. 
In
this case, one has that
$$
f_{_\cA} (\De)= \{\bar x^s\} \times [-1,1]^u \times \{\la\,x^c\}
 \quad \mbox{where} \quad \la\,
x^c\in \left(0, \frac{\mu}{\la-1}\right).
$$
Thus $f_{_\cA} (\De)$ is in between $W^{\st}_{\loc}(P,f)$ and
$W^{\st}_{\loc}(Q,f)$.

\medskip

\noindent{\bf Case 2:} $x^c\in I_2=(\frac{\mu}{\la},\frac{\mu}{\la-1})$.
Note that in this case one gets
$$
f_{_\cB}(\De)= \{\bar x^s\} \times [-1,1]^u \times \{\la\,x^c-\mu\}
\quad \mbox{where} \quad \la\, x^c-\mu \in \left(0,
\frac{\mu}{\la-1}\right).
$$
Hence $f_{_\cB} (\De)$ is in between $W^{\st}_{\loc}(P,f)$ and
$W^{\st}_{\loc}(Q,f)$.

This completes the proof of {\bf BH6)} for disks parallel to $E^\uut$.

\medskip

We next consider general $\uut$-disks. We begin with two claims.

\begin{clai}\label{c.1} Consider $\tau_1<\frac{\mu}{\la\, (\la-1)}$ and any $\uut$-disk $\De$
(i.e. tangent to the cone-field $\cC^\uut_\alpha$) through a point
$(x^s,x^u,x^c)$ with $x^c\leq \tau_1$. Then, for every
$\alpha\in(0,1)$ small enough,
\begin{itemize}
\item
  $f_{_\cA}(\De)$ is at the left
of $W^\st_\loc(Q)$,
\item
assume that $\De$ is at the right of $W^\st_\loc(P)$, then
$f_{_\cA}(\De)$ is in between $W^{\st}_{\loc}(P,f)$ and
$W^{\st}_{\loc}(Q,f)$.
\end{itemize}
\end{clai}
\begin{demo}
The first statement follows from the compactness of the set of
vertical disks with $x^c\leq \tau_1$ and the uniform
convergence in $\alpha$ of the $\uut$-disks to the vertical disks
(parallel to $E^\uut$) as $\alpha\to 0$. 

From the first part of the claim we know that $f_{_\cA}(\De)$ is at the left of 
$W^{\st}_{\loc}(Q,f)$. It remains to check that $f_{_\cA}(\De)$ is at the right of 
$W^{\st}_{\loc}(P,f)$. The disk $\De$ contains a point of the form $(x^s,0^u,x^c)$ with $x^c>0$.
Thus $f_{\cA}(\De)$ contains the point $(S_1(x^s),0^u,\la\,x^c)$. This implies that 
$f_{\cA}(\De)$ is at the right of 
$W^{\st}_{\loc}(P,f)$.
\end{demo}

Arguing as above and using Case 2, one also deduces the following:

\begin{clai}\label{c.2} Consider $\tau_2>
\frac{\mu}{\la}$
and any $\uut$-disk $\De$
 through a point
$(x^s,x^u,x^c)$ with $x^c\geq \tau_2$. Then, for every
$\alpha\in(0,1)$ small enough,
\begin{itemize}
\item
$f_{_\cB}(\De)$ is at the right of $W^\st_\loc(P)$,
\item
assume that $\De$ is at the left of $W^\st_\loc(Q)$, then
$f_{_\cB}( \De)$ is in between $W^{\st}_{\loc}(P,f)$ and
$W^{\st}_{\loc}(Q,f)$.
\end{itemize}
\end{clai}

To get condition {\bf BH6)},  note that, since that $(\la-1)\in(0,1)$, one has
$\frac{\mu}{\la}<\frac{\mu}{\la\, (\la-1)}$. Now it is enough to note that, for
$\alpha\in(0,1)$ small enough, for every point $(x^s,x^u,x^c)$ in
a $\uut$-disk $\De$ in between $W^{\st}_{\loc}(P,f)$ and
$W^{\st}_{\loc}(Q,f)$ one has:
\begin{itemize}
\item either $x^c< \frac{\mu}{\la\, (\la-1)}$ and then, by Claim~\ref{c.1}, $f_{_\cA}(\De)$ is in between $W^{\st}_{\loc}(P,f)$ and
$W^{\st}_{\loc}(Q,f)$,
\item or $x^c> \frac{\mu}{\la}$ and then, by Claim~\ref{c.2}, $f_{_\cB}
(\De)$ is in between $W^{\st}_{\loc}(P,f)$ and
$W^{\st}_{\loc}(Q,f)$.
\end{itemize}

We have checked  that the set $\La$ satisfies conditions {\bf
BH1)}--{\bf BH6)}. Therefore the set $\La$ is a blender-horseshoe
and the proof of Proposition~\ref{p.blenderhorseshoe} is complete.
\end{demo}

\subsection{Strong homoclinic intersections and generation of blender-horseshoes}\label{sss.affine}

In this section, we state the generation of blender-horseshoes by saddle-node and
flip periodic points with strong homoclinic intersections.

Let $f$ be a diffeomorphism  with a periodic point $S$ such that the
tangent bundle of $M$ at $S$ has a $Df^{\pi(S)}$-invariant dominated
splitting $T_S M=E^{\sst}\oplus_{_<} E^\ct \oplus_{_<} E^\uut$,
where $E^{\sst}$ is uniformly contracting, $E^\uut$ is uniformly
expanding, and $E^\ct$ is one-dimensional. We say that $S$ is a
{\emph{saddle-node}} (resp. {\emph{flip}}) if the eigenvalue of
$Df^{\pi(S)}$ corresponding to the (one-dimensional) central
direction $E^c$ is $1$ (resp. $-1$).

 Consider the
strong stable and unstable manifolds of the orbit of $S$ (denoted by $W^\sst(S)$ and $W^\uut(S)$).
We
say that the $S$ has a {\emph{strong homoclinic intersection}} if
there is a point $X\in W^\sst(S)\cap W^\uut(S)$ such that $X\ne
f^i(S)$ for all $i$. The point  $X$ is a {\emph{strong homoclinic
point of $f$}} associated to $S$.

%The construction in Section~\ref{ss.prototype} implies that for
%every $\lambda$ close to $1$, the map $f_{\lambda,0}$ has strong
%homoclinic intersections associated to the saddles
% $P=(0^s,0^u,0)\in \DD_1$ and $Q=(a^s,a^u,0)\in \DD_2$.

Let $f$ be a diffeomorphism with a strong homoclinic intersection
associated to a saddle-node $S$. In \cite[Section 4.1]{BDijmj} it
is shown that there are $k\ge 1$ and a
$C^1$-perturbation $g$ of $f$ such that $g^k$ has an  a affine
horseshoe $\Lambda$ associated to $S$ whose central expansion is
arbitrarily close to $1$ (recall Section~\ref{ss.prototype}).  
Considering perturbations similar
to the ones in Proposition~\ref{p.blenderhorseshoe}, \cite[Sections 4.1.1-2]{BDijmj} gives the following:

\begin{prop}\label{p.stronghomoclinic}
Consider a diffeomorphism $f$ having a strong homoclinic
intersection associated to a saddle-node  $S$. There is a
diffeomorphism $g$ arbitrarily $C^1$-close to $f$ with a
$\cut$-blender-horseshoe having $S$ as a reference saddle.

The same statement holds for $\cst$-blender-horseshoes, i.e.
$\cut$-blender-horseshoes for $f^{-1}$.
\end{prop}

We observe that in \cite{BDijmj} the terminology blender-horseshoe it is not used.
However, the constructions in \cite{BDijmj} provide prototypical blender-horseshoes exactly
as the ones in Section~\ref{ss.prototype}. In fact, these constructions are the main motivation
(and model) for our definition of blender-horseshoe. Thus Proposition~\ref{p.stronghomoclinic}
just reformulates some results in \cite{BDijmj} using this new terminology of blender-horseshoes.

For diffeomorphisms having flip points we need the following
lemma (see \cite[Remark 4.6]{BDijmj}):

\begin{lemm}\label{l.stronghomoclinic}
Consider a diffeomorphism $f$ having a strong homoclinic
intersection associated to a flip point $S$. There is a diffeomorphism
$g$ arbitrarily $C^1$-close to $f$
 having a saddle node $S'$ with a strong homoclinic intersection
and such that the orbit of $S'$ remains in an arbitrarily small
neighborhood of the orbit of the initial flip point $S$.
\end{lemm}
\begin{demo}
Consider a $2$-parameter family of deformations $f_{s,t}$ of
$f=f_{0,0}$ such that
\begin{itemize}
\item
the parameter $s$ corresponds to a (non-generic) unfolding of the
flip, generating a saddle-node $S'$ close to the flip $S$ of
period twice the period of $S$, and
\item
the parameter $t$ corresponds to the unfolding of the strong
homoclinic intersection of the flip $S$, the local strong stable
manifold of $S$ ``passing from the left to the right" of the local
unstable manifold of $S$.
\end{itemize}
Then for every $s\neq0$ small enough, there is a small parameter
$t=t(s)$, $t(s)\to 0$ as $s\to 0$, such that the saddle-node $S'$ has a strong homoclinic
intersection.
\end{demo}

\subsection{Co-index one cycles and blender
horseshoes}\label{sss.coindexblender}

In this section, we borrow some arguments and results from
\cite{ABCDW,BDijmj} in order to prove  that diffeomorphisms with co-index
one heterodimensional cycles yield  blender-horseshoes.

\begin{prop}\label{p.coindexblenders}
Let $f$ be a diffeomorphism with a heterodimensional
cycle associated to saddles $P$ and $Q$ with
$\inds(P)=\inds(Q)+1$.
%such that
%$s=\mbox{$\st$-index\,}(P)=\mbox{$\st$-index\,}(Q)+1$. Then
Then there is $g$ arbitrarily $C^1$-close to $f$  with a saddle
$R$ such that:
\begin{enumerate}
\item $\inds(R)=\inds(Q)$,
the orbit of $R$ has a dominated splitting
$E^\sst\oplus_{_<} E^\ct\oplus_{_<}E^\uut$ with three non-trivial
bundles such that  $E^\sst$ and $E^\uut$ are uniformly contracting
and expanding, respectively, $\dim (E^\sst)=\inds(Q)$, and $\dim (E^\ct)$ has dimension one and
is expanding,
\item
there is $\cut$-blender-horseshoe of $g$ having $R$ as a reference
saddle,
\item \label{i.intermediate} $W^\st(R)$  intersects transversely
$W^\ut (Q)$, and
\item \label{i.intermediate2}
$W^\uut(R)$ meets transversely $W^\st(P)$.
\end{enumerate}
\end{prop}

The arguments for proving this proposition  can be found scattered
along several constructions in \cite{BDijmj}. But, unfortunately, this result is not stated explicitly there and its prove involves some adaptations of the constructions in \cite{BDijmj}. 
As the proof of Proposition~\ref{p.coindexblenders} is somewhat technical, we next explain
the sequence of arguments we borrow from \cite{BDijmj} and their
adaptations in order to prove this proposition.
%On the other hand, the construction in \cite{BDijmj} does not take care (explicitly) about the intersections
% between
%$W^\st(R)$ and  $W^\ut (Q)$ and between $W^\uut(R)$ and
%$W^\st(P)$. However, a careful analysis of the the constructions in \cite{BDijmj} provides
%such intersections, obtaining  
%items~\ref{i.intermediate} and \ref{i.intermediate2} of
%Proposition~\ref{p.coindexblenders}. Finally, similar arguments to get
% these intersections can be found in \cite{ABCDW} (again
%without mentioning explicitly generation of blenders).
The proof of Proposition~\ref{p.coindexblenders} consists of several reductions to simpler cases
we proceed to explain. Let us
begin with a  definition.

\begin{defi}[Strong-intermediate saddles] Let $f$ be a diffeomorphism having two periodic saddles
$P$ and $Q$ of indices $\inds(P)=\inds(Q)+1$. A
periodic point  $R$  is \emph{strong-intermediate} with respect to
$P$ and $Q$, denoted by $Q\prec_{\ut,\sst} R\prec_{\uut,\st} P$,
if:
\begin{itemize}
\item the orbit of $R$ is partially hyperbolic and has a dominated splitting
$E^\sst\oplus_{_<} E^\ct\oplus_{_<}E^\uut$ with three non-trivial
bundles such that  $E^\sst$ and $E^\uut$ are uniformly contracting
and expanding, $\dim (E^\sst)=\inds(Q)$, and $\dim (E^\ct)=1$,
\item
the strong stable manifold of $R$ meets transversely the unstable
manifold of $Q$ and the strong unstable manifold of $R$ meets
transversely the stable manifold of $P$, in a formula,
$$W^\sst(R)\pitchfork W^\ut(Q)\neq \emptyset\quad\mbox{and}\quad W^\uut(R)\pitchfork W^\st(P)\neq
\emptyset.$$
\end{itemize}
\end{defi}

 Note that if a (hyperbolic) saddle $R$ with $\inds (R)=\inds (Q)$ is strong-intermediate
 with respect  to $P$ and $Q$ then it
satisfies items~\ref{i.intermediate} and \ref{i.intermediate2} in
Proposition~\ref{p.coindexblenders}. 

We need the following lemma:

\begin{lemm}\label{l.strongintermediate}
Consider two saddles  $P$ and $Q$   in the same chain recurrence class
$C$ and
 a periodic point $R$ which is
strong-intermediate with respect to $P$ and $Q$.  Then $R\in C$.
\end{lemm}

\begin{demo}
We construct a pseudo-orbit going from $R$ to $P$
(the other pseudo-orbits are obtained similarly).
Take a point
$X\in W^\sst(R)\cap W^\ut(Q)$. Note that there are arbitrarily large $n$ and $m$ such that
$\{f^{-n}(X),\dots, f^{m}(X)\}$ is a segment of orbit
starting (arbitrarily) close to $R$ and ending close to $Q$. Since $Q$ and
$P$ are in the same chain recurrent class, there is a finite
pseudo-orbit going from $Q$ to $P$. A pseudo-orbit going
from $R$ to $P$ is obtained concatenating these two
pseudo-orbits. This concludes the sketch of the proof of the lemma.
\end{demo}

We now explain the generation of strong-intermediate saddles.

\subsubsection{Reduction to the case  of  cycles associated to saddles with real central eigenvalues}
Given a periodic point $R$ of a diffeomorphism $f$, write
$\la_1(R),\dots,\la_n(R)$ the eigenvalues of $Df^{\pi(R)} (R)$
counted with multiplicity and ordered in increasing modulus
($|\la_i(R)|\le |\la_{i+1}(R)|$). We say that $\la_i(R)$ is the
$i$-th multiplier of $R$.

Consider a diffeomorphism $f$ having a co-index one  cycle
associated to period points $A$ and $B$ with
$\inds(P)=\inds(Q)+1=s+1$. The cycle \emph{has real central
eigenvalues} if $\la_{s+1}(A)$ and $\la_{s+1}(B)$ are both real
and
$$|\la_s(A)|<|\la_{s+1}(A)|<1<|\la_{s+2}(A)|\quad\mbox{and}\quad
|\la_s(B)|<1<|\la_{s+1}(B)|<|\la_{s+2}(B)|.$$

Before proving Proposition~\ref{p.coindexblenders},
we recall the
following two facts: 

\begin{description}
\item{\bf Fact 1:}
 \cite[Theorem 2.1]{BDijmj} claims that, if $f$ has a co-index
one cycle associated to $A$ and $B$ then there is $g$ arbitrarily
$C^1$-close to $f$ with a co-index one cycle with real central
eigenvalues. Moreover, this cycle can be chosen  associated to
saddles $A'_g$ and $B'_g$ homoclinically related to the
continuations $A_g$ and $B_g$ of $A$ and $B$, respectively.
\item{\bf Fact 2:}
Assume  that there is  a diffeomorphism $h$ arbitrarily
$C^1$-close to $g$ having a $\cut$-blender-horseshoe $\La$ with a
reference saddle $R_h$ which is strong-intermediate to $A'_h$ and $B'_h$.
Since two saddles being homoclinically related is a $C^1$-robust
relation, one has $R_h$ is strong-intermediate with respect to
$A_h$ and $B_h$. In this case, the proof of
Proposition~\ref{p.coindexblenders} is complete.
\end{description}

In view of these two facts, to prove Proposition~\ref{p.coindexblenders} it is
enough to consider the case where the saddles $P$ and $Q$  in the cycle have 
real central eigenvalues and to
check that these cycles generate  strong-intermediate saddles as
in Fact 2, see Proposition~\ref{p.saddleblenders}. We now go to the details of this construction.

\subsubsection{Reduction to the generation of  saddle-node or flip points}
In this section, we show that
Proposition~\ref{p.coindexblenders} is a consequence of the
following result.

\begin{prop}\label{p.saddleblenders}
Let $f$ be a diffeomorphism with a co-index one cycle associated to
saddles $P$ and $Q$ with real central eigenvalues.
%such that
%$s=\mbox{$\st$-index\,}(P)=\mbox{$\st$-index\,}(Q)+1$. Then
Then there is $g$ arbitrarily $C^1$-close to $f$  having a saddle-node
or flip periodic point $R_g$ such that:
\begin{itemize}
\item $R_g$ has a strong homoclinic intersection and
\item  $R_g$ is strong-intermediate to $P_g$ and $Q_g$.
\end{itemize}
\end{prop}

This proposition  is a stronger version of \cite[Theorem 2.3]{BDijmj}, 
adding the strong-intermediate property of $R$ with respect to $P$
and $Q$.

%We now deduce 
%Proposition~\ref{p.coindexblenders} from
%Proposition~\ref{p.saddleblenders}.  

\begin{demo}[Proposition~\ref{p.saddleblenders} implies Proposition~\ref{p.coindexblenders}]
First, observe that the strong
unstable and strong stable manifolds of $R_g$ depend continuously on the
diffeomorphism (while there is defined a continuation of $R_g$).

If $R=R_g$ is a saddle-node then
Proposition~\ref{p.stronghomoclinic} gives a diffeomorphism $h$
arbitrarily $C^1$-close to $g$ (thus arbitrarily close to $f$)
with a blender-horseshoe having $R$ as a reference saddle. Therefore, for $h$
close enough to $g$, one gets the announced intersections between
the strong invariant manifolds (intermediate intersections).

If $R$ is a flip then Lemma~\ref{l.stronghomoclinic} gives a
perturbation $h$ of $g$ with a saddle-node with a strong
homoclinic intersection and with the strong-intermediate property.
Thus we are in the first case. This completes the proof of our claim.
\end{demo}

Therefore it is enough to prove Proposition~\ref{p.saddleblenders}. The proof of 
the proposition  is similar to the one of  \cite[Theorem 2.3]{BDijmj} and
consists of
several steps. We next explain and adapt these steps.

\subsubsection{Reduction to the creation of weak hyperbolic saddles}
We now see that Proposition~\ref{p.saddleblenders} (hence of
Proposition~\ref{p.coindexblenders}) follows from:

\begin{prop}\label{p.weak} Let $f$ be a diffeomorphism having a
co-index one cycle associated to saddles  $P$ and $Q$, $\inds
(P)=\inds(Q)+1$, with real central eigenvalues. Then there are a
constant $C>1$ and a sequence of diffeomorphisms $(f_n)$, $f_n
\stackrel{C^1}{\longrightarrow} f$, such that every $f_n$ has a
periodic point $R_n$ such that:
\begin{itemize}
\item $R_n$ has a one-dimensional center-unstable direction whose
corresponding multiplier $\la_n^c$ satisfies
$|\la_n^c|\in[\frac1C,C]$,
\item
$W^\uut(R_n)$ and $W^\sst(R_n)$ have a quasi-transverse
intersection, therefore $R_n$ has strong homoclinic  intersections
\item  the periods of $R_n$ go to infinity as $n\to \infty$, and
\item
$R_n$ is strong-intermediate with respect to $P_n$ and $Q_n$ (the
continuations of $P$ and $Q$ for $f_n$).
\end{itemize}
\end{prop}

This proposition  is a stronger version \cite[Proposition
3.3]{BDijmj}, adding the intersection property of the strong
invariant manifolds.

\begin{demo}[Proposition~\ref{p.weak} implies Proposition~\ref{p.saddleblenders}]
We proceed exactly as in \cite[page
484]{BDijmj} (proof of Theorem 2.3 using Proposition 3.3). 
We just
perform a local  $C^1$-perturbation of $f_n$ supported in a
small neighborhood of $R_n$, turning the the central
eigenvalue of $R_n$ equal to $\pm 1$ while keeping the strong
homoclinic point of $R_n$ and the transverse intersections
$W^\sst(R_n)\pitchfork W^\ut (Q_n)\neq \emptyset$ and
$W^\uut(R_n)\pitchfork W^\st(P_n)\neq\emptyset$. 

In this way, we
get diffeomorphisms $g$ (arbitrarily close to $f$) with
saddle-node or flip points $R_g$ (depending if $\la_n^c$ is
positive or negative) with  strong homoclinic intersections and
being strong-intermediate with respect to $P_g$ and $Q_g$. This
completes the proof of Proposition~\ref{p.saddleblenders}.
\end{demo}

\subsection{Proof of Proposition~\ref{p.weak}}
\label{sss.profofweak}
The following steps of the proof of \cite[Proposition 3.3]{BDijmj}
are described in \cite[page 484]{BDijmj}:

\begin{description}
\item{\bf Step 1:}
One first puts the cycle in a kind of normal
form called \emph{simple cycle}. In fact, \cite[Proposition
3.5]{BDijmj} implies that, after a
$C^1$-perturbation, one can assume  that the cycle is simple.
\item{\bf Step 2:}
One shows that the dynamics in a simple cycle is given
(up to a renormalization) by a model family, denoted by
$F^{\pm,\pm}_{\la,\be,t}$. Moreover,  perturbations of this model
family correspond to perturbations of the initial dynamics.
\end{description}

Therefore, to prove Proposition~\ref{p.weak}, it is enough to
consider model families $F^{\pm,\pm}_{\la,\be,t}$ and their
perturbations. Hence it is enough to adapt the perturbations of
these normal families  in order to get the intersection
properties  between the strong invariant manifolds.

Proposition 3.8 in
\cite{BDijmj} claims that the unfolding of  co-index one cycles generates sequences of saddles $A_{n,m}$
whose orbits are contained in a neighborhood of the cycle. We now see that these saddles can be taken
with the strong-intermediate property (relative to saddles in the initial cycle).

\begin{lemm}\label{l.intermediate} In {\em \cite[Proposition 3.8]{BDijmj}}, for every integer
$n,m$ large enough, (in fact larger that the integer $N$ in the
statement), all  periodic  points $A_{n,m}$ in the proposition are
strong-intermediate with respect to $P$ and $Q$.
\end{lemm}
Note that we have the following string of implications:
$$
\mbox{Lemma~\ref{l.intermediate}}\Rightarrow
\mbox{Proposition~\ref{p.weak}}\Rightarrow
\mbox{Proposition~\ref{p.saddleblenders}} \Rightarrow
\mbox{Proposition~\ref{p.coindexblenders}}.
$$
Therefore to prove  Proposition~\ref{p.coindexblenders} it remains to prove Lemma~\ref{l.intermediate}.

\begin{demo} The model maps $F^{\pm,\pm}_{\la,\be,t}$ are 
defined  on some cubes and  their restrictions to each of
these cubes are affine maps $\cA_\lambda$, $\cB_\beta$,
$\cT^\pm_{1,t}$, and $\cT^\pm_2$
which preserve a constant dominated splitting: the strong stable bundle
is the horizontal space $\RR^s\times \{(0,0^u)\}$, the strong unstable bundles  is the vertical one
$ \{(0^s,0)\}\times \RR^u$,
and the center bundle is one dimensional $\{0^s\}\times \RR \times \{0^u\}$. More precisely 
(see Figure~\ref{f.trans}):

\begin{figure}[htb]

\begin{center}
\psfrag{Yp}{$Y$}
\psfrag{Xp}{$f^{N_2}(X)$}
\psfrag{Yq}{$f^{N_1}(Y)$}
\psfrag{Xq}{$X$}
\psfrag{Yq}{$f^{N_1}(Y)$}
\psfrag{Xq}{$X$}

\psfrag{A}{$\cA_\lambda$}

\psfrag{B}{$\cB_\beta$}

\psfrag{Es}{$E^{ss}$}
\psfrag{Eu}{$E^{uu}$}
\psfrag{Ec}{$E^{c}$}
%\psfrag{wsp}{$W^s(P)$}
%\psfrag{wup}{$W^u(P)$}
%\psfrag{wsq}{$W^s(Q)$}
%\psfrag{wuq}{$W^u(Q)$}
\psfrag{P}{$P$}
\psfrag{Q}{$Q$}
\psfrag{T1}{$\mfT_{1,t}$}
\psfrag{T2}{$\mfT_2$}
\includegraphics[height=1.8in]{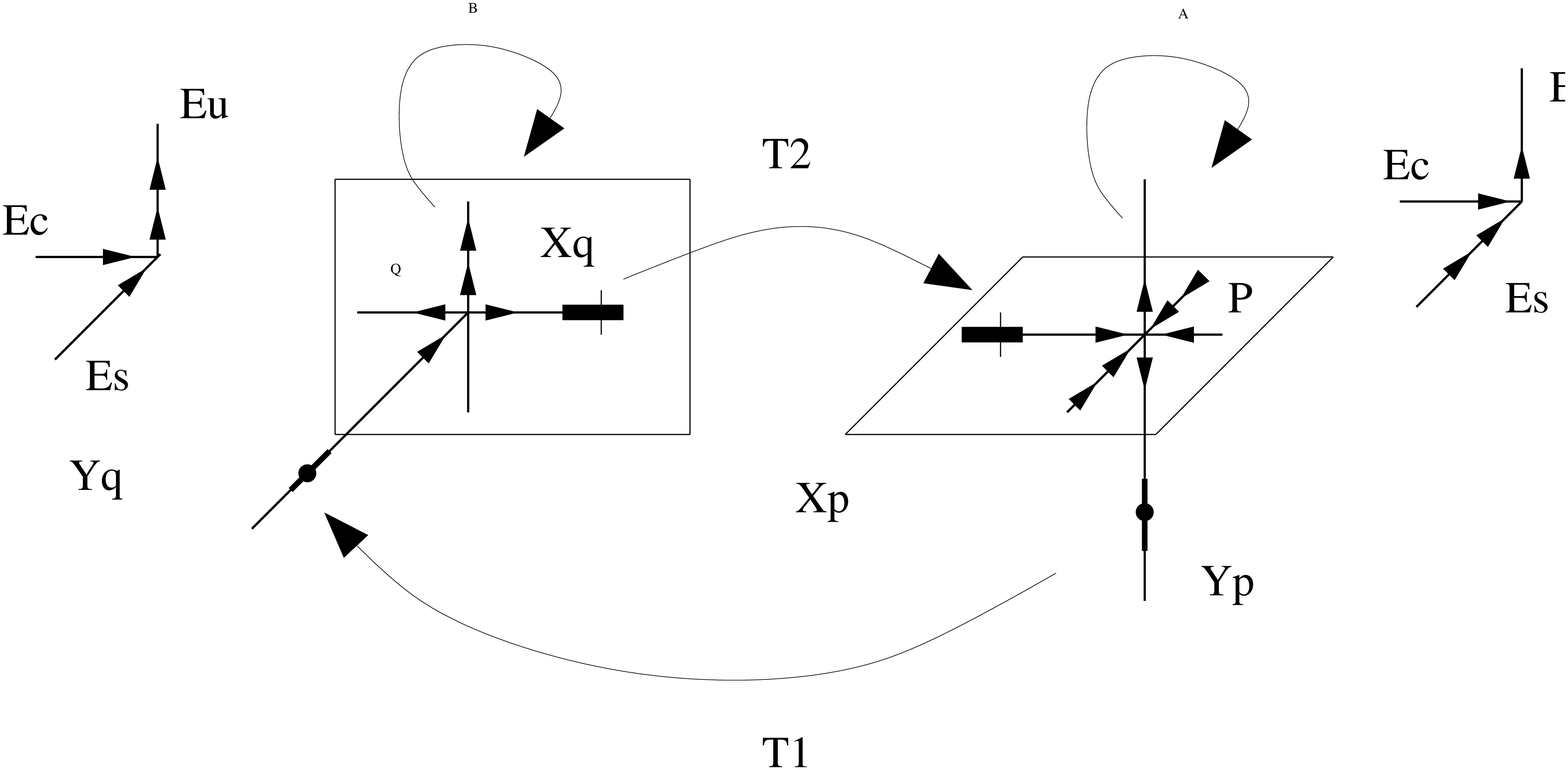}
 \caption{The model maps $F^{\pm,\pm}_{\la,\be,t}$.}
\label{f.trans}
\end{center}

\end{figure}

\begin{itemize}
 \item 
The
maps $\cA_\lambda$ and $\cB_\be$ are  $Df^{\pi(P)}(P)$ and
$Df^{\pi(Q)}(Q)$ and correspond to iterates of $f$ close to $P$
and $Q$, respectively. 
\item
The sub-scripts $\la$ and $\be$ correspond
to the central multipliers of $Df^{\pi(P)}(P)$ and
$Df^{\pi(Q)}(Q)$. 
\item
The maps $\cT^\pm_{1,t}$ and $\cT^\pm_2$ are the
{\emph{transitions}} of the cycle. The map $\cT^\pm_2$ corresponds to a
fixed number $N_2$ of iterates from a neighborhood of $Q$ to  a
neighborhood of $P$ following a segment of orbit of a transverse heteroclinic
point $X$ in $W^\st(P)\pitchfork W^\ut(Q)$. 

Similarly, the map $\cT^\pm_{1,t}$
corresponds to a fixed number $N_1$ of iterates from  a neighborhood of
$P$ to a neighborhood of $Q$ following a segment of orbit of a fixed quasi-transverse
heteroclinic point $Y$ in $W^\ut(P)\pitchfork W^\st(Q)$. 
The parameter $t$ of $\cT^\pm_{1,t}$ corresponds to the unfolding of the cycle.
\item
The
super-script $\pm$ is positive if the transition map preserves the
orientation in the central bundle and negative if otherwise.
\end{itemize}
For details see
\cite[page 488]{BDijmj}.

\begin{figure}[htb]

\begin{center}
\psfrag{Yp}{$Y$}
\psfrag{Xp}{$f^{N_2}(X)$}
\psfrag{Yq}{$f^{N_1}(Y)$}
\psfrag{Xq}{$X$}
\psfrag{Yq}{$f^{N_1}(Y)$}
\psfrag{Xq}{$X$}

\psfrag{A}{$\cA_\lambda$}

\psfrag{amn}{$A_{m,n}$}

\psfrag{B}{$\cB_\beta$}

\psfrag{Es}{$E^{ss}$}
\psfrag{Eu}{$E^{uu}$}
\psfrag{Ec}{$E^{c}$}
\psfrag{wsp}{$W^s(P)$}
\psfrag{wuq}{$W^u(Q)$}
%\psfrag{wsq}{$W^s(Q)$}
%\psfrag{wuq}{$W^u(Q)$}
\psfrag{P}{$P$}
\psfrag{Q}{$Q$}
\psfrag{T1}{$\mfT_{1,t}$}
\psfrag{T2}{$\mfT_2$}
\includegraphics[height=2.0in]{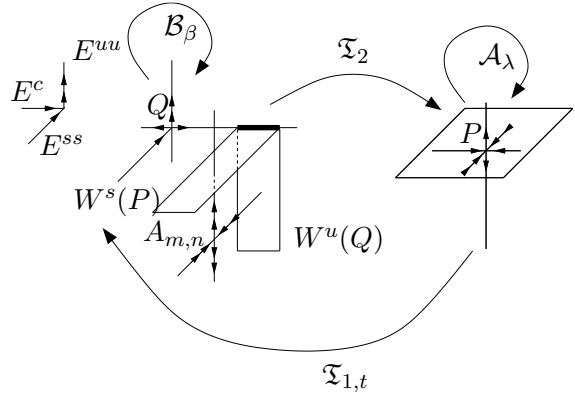}
 \caption{The intermediate saddle $A_{m,n}$.}
\label{f.amn}
\end{center}

\end{figure}

By definition,  the points $A_{m,n}=(a^s,a,a^u)$ are fixed points of the
composition $\cB_\beta^n\circ \cT^\pm_{1,t}\circ
\cA_\lambda^m\circ\cT^\pm_2$. In particular, the point $A_{m,n}$
belongs to domain of definition $\Si_Q$  of $\cT^\pm_2$, see Figure~\ref{f.amn}. This domain
(defined in \cite[page 488]{BDijmj}) is the cube
$$
\Si_Q=[-1,1]^s\times [b_Q-\delta,b_Q+\delta]\times[-1,1]^u,
$$
where
$$
[-1,1]^s\times [b_Q-\delta,b_Q+\delta]\times\{0^u\} \subset
W^\st(P)\quad \mbox{and} \quad \{0^s\}\times
[b_Q-\delta,b_Q+\delta]\times[-1,1]^u\subset W^\ut(Q).
$$
Moreover, the local strong stable and strong unstable manifolds of
the point $A_{m,n}=(a^s,a,a^u)$ are (see \cite[page 495, first
paragraph]{BDijmj}):
$$ W^\sst_\loc(A_{m,n})= [-1,1]^s\times \{(a,a^u)\}
\quad\mbox{and}\quad W^\uut_\loc(A_{m,n})= \{(a^s,a)\}\times
[-1,1]^u.
$$
This means  that 
$$
 W^\sst(A_{m,n}) \pitchfork W^\ut(Q) \ne\emptyset \quad \mbox{and}
\quad
 W^\uut_\loc(A_{m,n}) \pitchfork W^\st(P) \ne\emptyset.
$$
Therefore
the saddle
$A_{m,n}$ is strong-intermediate with respect to
$P$ and $Q$, ending the proof of the lemma.
\end{demo}

\section{Robust tangencies and heterodimensional cycles for $C^1$-generic diffeomorphisms}
\label{s.tmain}

In this section, we prove Theorem~\ref{t.main}. First, in
Section~\ref{ss.generic}, we state some properties about
$C^1$-generic diffeomorphisms. In Section~\ref{ss.generation}, we state
the $C^1$-generic occurrence of blender-horseshoes in homoclinic classes
containing saddles of different indices
(Theorem~\ref{th.genericblender}). Finally, in
Section~\ref{s.homoclinic}, we state the existence of robust
homoclinic tangencies inside homoclinic classes with index
variation and lack of domination (Proposition~\ref{p.main}), 
completing the proof of Theorem~\ref{t.main}. We close this paper presenting and
extension of   \cite[Theorem 1.16]{BDijmj}
about the occorrence of robust heterodimensional cycles inside
chain recurrence classes (see Theorem~\ref{t.heterocycle} in Section~\ref{ss.complemet}).

\subsection{$C^1$-generic properties of $C^1$-diffeomorphisms}
\label{ss.generic} We now collect some properties of
$C^1$-generic diffeomorphisms.
 According to \cite[Remarque 1.10]{BoCroIn} and
\cite[Theorem 1]{ABCDW}, there is a residual subset set $\cG$ of
$\diffM$ such that, for every $f\in\cG$,
\begin{itemize}
\item every periodic point of $f$ is hyperbolic,
\item for every periodic point $P$ of $f$,  its homoclinic class $H(P,f)$
and its chain recurrence class $C(P,f)$ are equal,
\item any homoclinic class $H(P,f)$ containing periodic points of
$\ut$-indices $\alpha$ and $\beta$ also contains saddles of
$\ut$-index $\tau$, for every $\tau\in[\alpha,\beta]\cap\NN$.
\end{itemize}

\begin{lemm}{\bf (\cite[Lemma 2.1]{ABCDW})}\label{l.ABCDW} There is a residual subset $\cG_0\subset \cG$ of $\diffM$ such
that,
for every $f\in\cG_0$ and every pair of periodic points $P_f$ and
$Q_f$ of $f$, there is a neighborhood $\cU_f$ of $f$ in $\cG_0$
such that:
\begin{itemize}
\item either $H(P_g,g)=H(Q_g,g)$ for all $g\in\cU_f\cap\cG_0$,
\item or $H(P_g,g)\cap H(Q_g,g)=\emptyset $ for all $g\in\cU_f\cap \cG_0$.
\end{itemize}
\end{lemm}

\begin{rema}[Proof of Claim 2.2 in \cite{ABCDW}]
Using a
filtration given by Conley theory, one has that the property of
two hyperbolic saddles to be in different chain recurrence classes
is $C^1$-robust.
\label{r.conley}
\end{rema} 

Next lemma claims that, for $C^1$-generic
diffeomorphisms, the property of the chain recurrence classes of
two periodic points  to be equal is also a $C^1$-robust property.

\begin{lemm}\label{l.robustclass} Let $\cG_0$ be the residual set of of $\diffM$
in Lemma~\ref{l.ABCDW}.
For every $f\in\cG_0$ and every pair of periodic  points $P_f$ and
$Q_f$ of $f$, the property of $Q_f$ belonging to the chain
recurrence class $C(P_f,f)$ of $P_f$ is $C^1$-robust: 
if $Q_f\in C(P_f,f)$ then
$Q_g\in C(P_g,g)$ for all
$g$ $C^1$-close to $f$.
\end{lemm}
\begin{demo} Let $f\in\cG_0$ and suppose that $Q_f\in C(P_f,f)$. Since
$f\in\cG$ one has that $C(P_f,f)=H(P_f,f)$. As $f\in\cG_0$, by
Lemma~\ref{l.ABCDW}, there is a $C^1$-open neighborhood $\cU_f$ of
$f$ such that $H(P_g,g)=H(Q_g,g)$ for every $g\in\cU_f\cap\cG_0$. In
particular, $Q_g\in H(P_g,g)$ for all $g\in\cU_f\cap\cG_0$.

Assume that there is $g\in\cU_f$ such that $Q_g\notin C(P_g,g)$. By
Remark~\ref{r.conley}, one has that $Q_h\notin C(P_h,h)$ for
every
 $h$ in a small
neighborhood $\cV_g$ of $g$ contained in $\cU_f$. Choosing
$h\in\cG_0\cap\cV_g$ one gets that $Q_h\not\in C(P_h,h)=H(P_h,h)$,
 a contradiction.
\end{demo}

\subsection{Generation of blender-horseshoes}\label{ss.generation}

In this section, we prove that blender-horseshoes occur
$C^1$-generically for homoclinic classes with index variability (i.e., containing saddles with
different indices).

\begin{theor}
\label{th.genericblender} There is a residual subset $\cR$ of
$\diffM$  of diffeomorphisms $f$ such that  for every
homoclinic class $H(P,f)$  containing a
hyperbolic  saddle $Q$ with $\inds(Q)>\inds(P)$ there is a
transitive hyperbolic set $\Si$ containing $P$ and a
$\cut$-blender-horseshoe $\La$. 
\end{theor}

Applying Theorem~\ref{th.genericblender} to
$f^{-1}$ one gets the following.

\begin{rema}\label{r.genericblender} Under the assumptions of
Theorem~\ref{th.genericblender}, every $C^1$-generic
diffeomorphism $f$ has a $\cst$-blender-horseshoe containing $Q$
and contained in a transitive hyperbolic set.
\end{rema}

We first prove a version of
Theorem~\ref{th.genericblender} for a  given fixed saddle
$P_f$:

\begin{prop}\label{pr.genericblender} Let $\cU$ be an open subset of $\diffM$ and $f\mapsto
P_f$ be a continuous map defined on $\cU$ associating to each
$f\in \cU$ a hyperbolic periodic point $P_f$ of $f$.

There is a residual subset $\cR=\cR_{P}$ of $\diffM$ with the
following property. For every diffeomorphism $f\in \cR\cap\cU$
such that $H(P_f,f)$ contains a saddle $B_f$ with
$\inds(B_f)>\inds(P_f)$ there is a transitive hyperbolic set
$\Si_f$ containing $P_f$ and a $\cut$-blender-horseshoe $\La_f$.
\end{prop}

Theorem~\ref{th.genericblender} will follow from this proposition
using  standard genericity arguments (the details can be found in the
end of this subsection).

\begin{demo}[Proof of Proposition~\ref{th.genericblender}]
Consider the residual subset $\cG_0$ of $\diffM$ in Lemma~\ref{l.ABCDW}.
Given any $f\in\cG_0\cap\cU$,  if $H(P,f)$ contains a  saddle $B_f$
with $\inds(B_f)>\inds(P_f)$, then there is a saddle
$Q_f\in H(P_f,f)$ with $\inds(Q_f)=\inds(P_f)+1$. Furthermore,
by Lemma~\ref{l.robustclass}, if $f\in\cG_0$ the point
$Q_f$  belongs robustly to the chain recurrence class $C(P_f,f)$. Moreover,  $C(P_f,f)=H(P_f,f)$.

 Let $\cW$ be  the set of diffeomorphisms:
$$\cW=\left\{ f\in \cU \colon \mbox{there is a saddle  $Q_f\in\Perf$
with} \quad
\begin{array}{ll}
&\inds(Q_f)=\inds(P_f)+1,\\
&\mbox{and}\\
&\mbox{$Q_f$ is $C^1$-robustly in $C(P_f,f)$}
\end{array}
\right\}
$$

By definition, the set $\cW$ is an open subset of $\cU$. Let
$$
\cU_1=\cU\setminus\overline{\cW}.
$$ 
By construction the set $\cU_1$ is open and 
$\cU_1\cup\cW$ is dense in $\cU$.

\begin{clai}  Let $f\in\cG_0\cap \cU_1$. Then $C(P_f,f)$ does not contain
any hyperbolic periodic point $\tilde Q$ with $\inds(\tilde
Q)>\inds(P_f)$.
\end{clai}
\begin{demo} The proof is by contradiction. Suppose that there is a saddle $\tilde Q\in C(P_f,f)$
with $\inds(\tilde Q)>\inds(P_f)$. As $f\in\cG_0$, there is a saddle
$Q_f\in C(P_f,f)=H(P_f,f)$ with $\inds(Q_f)=\inds(P_f)+1$. As
$f\in\cG_0$, the saddle $Q_f$ belongs $C^1$-robustly to $C(P_f,f)$.
Therefore $f\in \cW$, contradicting the definition of $\cU_1$.
\end{demo}

Thus,  by Lemma~\ref{l.robustclass}, it is enough to prove Proposition~\ref{pr.genericblender}  
for diffeomorphisms in $\cW$.
In other words, next lemma implies Proposition~\ref{pr.genericblender}.

\begin{lemm} The open set of diffeomorphisms $f$ having a transitive
hyperbolic set $\Si_f$ containing $P_f$ and a
$\cut$-blender-horseshoe $\La_f$ is dense in $\cW$.
\end{lemm}
\begin{demo} Consider $f\in\cW$ and a saddle $Q_f$ with
$\inds(Q_f)=\inds(P_f)+1$ which belongs robustly to $C(P_f,f)$.
Next lemma is an immediate
consequence of 
Hayashi connecting lemma in \cite{Ha}:
%it can be found in \cite{ABCDW}. Here we formulate the lemma
%in \cite{ABCDW} for chain recurrence classes instead of
%homoclinic classes, but the proof is the same. 
\begin{lemm}
[\cite{Ha}] 
Let $f$ be a diffeomorphism having
a pair of (hyperbolic) saddles $A_f$ and $B_f$ whose orbits are different and such that
$B_f\in H(A_f,f)$. Then there is $g$
arbitrarily $C^1$-close to $f$ such that $W^\ut(A_g,g)\cap W^\st(B_g,g)\ne\emptyset$.
\label{l.hayashihomoclinic}
\end{lemm}
Note that as $f\in \cW$ we have $Q_f\in C(P_f,f)$ in a robust way. After a first perturbation, we
can assume that $f\in \cG_0$ so that $Q_f\in H(P_f,f)$.
We
apply Lemma~\ref{l.hayashihomoclinic}  to the saddles $P_f$ and $Q_f$
to get a diffeomorphism $g$ close to $f$ with
a transverse intersection 
between  $W^\ut(P_g,g)$ and $W^\st(Q_g,g)$. Note that this transverse intersection persists after perturbation, so 
that we can assume that $g\in \cW \cap \cG_0$.
After a new application of Lemma~\ref{l.hayashihomoclinic} we can suppose that 
 $W^\st(P_g,g)\cap W^\ut(Q_g,g)\ne \emptyset$.
Therefore there is $g$
arbitrarily $C^1$-close to $f$, having a co-index
one  cycle associated to $P_g$ and $Q_g$. 

Applying
Proposition~\ref{p.coindexblenders} to the diffeomorphism $g$ with a co-index one cycle (associated to $P_g$ and $Q_g$), we get 
$h$ 
close to $g$, thus close to $f$, with a periodic point $R_h$ such that
\begin{itemize}
 \item 
$R_h$
is strong-intermediate with respect to $P_h$ and $Q_h$,
\item
$R_h$
has
the same index as $P_h$, and
\item
$R_h$
 is a reference saddle of a
$\cut$-blender-horseshoe $\La_h$.
\end{itemize}
Recall that, by
Lemma~\ref{l.continuation}, the continuation of this
blender-horseshoe is defined in a neighborhood of $h$. Thus
we can assume that $h\in \cG_0$.

As the saddle $R_h$ is (robustly) intermediate with respect  to
$P_h$ and $Q_h$ and $Q_h\in C(P_h,h)$, 
Lemma~\ref{l.strongintermediate} implies  that $R_h\in C(P_h,h)$.
Thus, from Lemma~\ref{l.robustclass}, we have that  $R_h$ belongs
robustly to $C(P_h,h)$. As $R_h$ and $P_h$ have the same index, Lemma~\ref{l.hayashihomoclinic}
gives a perturbation $\vfi$ of $h$  with transverse
intersections between the invariant
manifolds of these saddles.  That is, the saddles
 $R_\vfi$ and $P_\vfi$ are homoclinically
related, thus their homoclinic classes are equal. Therefore there is a transitive hyperbolic set
$\Sigma_\vfi$ containing $P_\vfi$ and the $\cut$-blender-horseshoe
$\La_\vfi$ (the continuation of $\La_h$). This ends the proof of
the lemma.
\end{demo}

The proof of Proposition~\ref{pr.genericblender} is now complete.
\end{demo}

\begin{demo}[Proof of Theorem~\ref{th.genericblender}]
Given a diffeomorphism $f$, denote by $\Perfn$ the set of periodic
points $P$ of $f$ of period $\pi(P)\le n$. To prove
Theorem~\ref{th.genericblender} it is enough to see that, for
every $n\in\NN$, there is a residual set $\cR_{\leq n}$ of
diffeomorphisms $f$ such that the conclusion of the theorem
holds for every periodic orbit $P\in \Perfn$. Then it is enough
to take $\cR=\cap_n \cR_{\le n}$.

 Note that there is a
$C^1$-open an dense subset $\cO_n\subset \diffM$ of
diffeomorphisms $f$ such that every periodic point $P\in \Perfn$
is hyperbolic. In particular, the cardinal of $\Perfn$  is finite
an locally constant in $\cO_n$. Moreover, each periodic point
$P\in \Perfn$ has a hyperbolic continuation in each (open)
connected component $\cU$ of the open set $\cO_n$. That is, there
are a constant $k=k(\cU)$ and continuous maps $f\mapsto P_{i,f}$,
$i=1,\dots, k$, such that $\Perfn=\{P_{1,f},\dots,P_{k,f}\}$, for
every $f\in\cU$.

Note that the set $\mbox{cc}(\cO_n)$ of connected components of
$\cO_n$ is countable. Therefore to prove
Theorem~\ref{th.genericblender} for periods $\pi \leq n$
 it is enough to see that  this result holds  in each connected component
$\cU$ of $\cO_n$.
  More precisely, for each connected component $\cU$ of $\cO_n$, we
first build a residual subset $\widetilde{\cR_\cU}$ of $\diffM$
such that the conclusion holds in the set $\widetilde{\cR_\cU}\cap
\cU$. We  now consider the set
$$
\cR_\cU=\widetilde{\cR_\cU}\cup (\cO_n\setminus \cU).
$$ Note
that the set $\cO_n\setminus \cU$ is the union of the open
connected components of $\cO_n$ different from $\cU$. Thus the set
$\cO_n\setminus \cU$ is open (and closed) in $\cO_n$ and therefore the set
$\cR_\cU$ is residual in $\diffM$. Finally, the announced residual
subset  $\cR_{\leq n}$ of $\diffM$ is the countable intersection
of the residual subsets $\cR_\cU$ of $\diffM$:
$$\cR_{\leq n}=\bigcap_{\cU\in \mbox{cc}(\cO_n)}\cR_\cU.$$

To complete the proof of the theorem it remains to define $\widetilde{\cR_\cU}$
for each component $\cU$ of $\cO_n$.
 Given $f\in\cU$ consider $\Perfn=\{P_{1,f},\dots,P_{k,f}\}$. For
each $i=1,\dots, k$, Proposition~\ref{pr.genericblender} gives a
residual subset $\cR_{P_i}$ of $\diffM$ where the conclusion
holds. The residual set $\widetilde{\cR_\cU}$ is the
finite intersection of the residual sets $\cR_{P_i}$. The proof of
Theorem~\ref{th.genericblender} is now complete.
\end{demo}

\subsection{Robust homoclinic tangencies under lack of domination}
 \label{s.homoclinic}

In this section, we conclude the proof of Theorem~\ref{t.main}
about $C^1$-generic existence of robust homoclinic tangencies inside
homoclinic classes with index variation and lack of domination. 
We first recall a key result stating the relation between lack of
domination and homoclinic tangencies.

\begin{theor}[Theorem 1.1 in \cite{gourmelon}]
Let  $P_f$ be a saddle of a diffeomorphism $f$ such that
the stable/unstable splitting defined over the set of
periodic points homoclinically related with $P_f$ is not
dominated. Then  there is
a diffeomorphism $h$ arbitrarily $C^1$-close to $f$ with a
homoclinic tangency associated to $P_h$.
\label{t.gourmelon}
\end{theor}

As
in Section~\ref{ss.generation}, we begin with a version of
Theorem~\ref{t.main} for a given fixed saddle.

\begin{prop}\label{p.main} Consider a diffeomorphism $g$ and
a hyperbolic saddle $P_g$ of $g$.
Let $\cU$ be an open subset of $\diffM$ such that the map
$f\mapsto P_f$ ($P_f$ a hyperbolic saddle) is continuous and well
defined. Then there is a residual subset $\cG_\cU$ of $\cU$ with
the following property. Let $f\in\cG_\cU$ be any diffeomorphism such
that:
\begin{itemize}
\item
the chain recurrence class $C(P_f,f)$ has a periodic point $Q_f$
with $\inds(Q_f)>\inds(P_f)$ and
\item
$C(P_f,f)$ does not admit a dominated splitting $E\oplus_{_<}F$ with
$\dim (E)=\inds(P_f)$.
\end{itemize}
The $C(P_f,f)$ has a
 transitive hyperbolic set containing
$P_f$ with a $C^1$-robust homoclinic tangency.
\end{prop}

\begin{demo}
Let $\cW_0\subset \cU$ be the set of diffeomorphisms $f$ such that
the chain recurrence class $C(P_f,f)$ of $P_f$ contains robustly a hyperbolic
periodic point $Q_f$ with $\inds(Q_f)>\inds(P_f)$ (i.e.,
$Q_g\in C(P_g,g)$ for all $g$ $C^1$-close to $f$).  By
definition and Remark~\ref{r.conley}, the set $\cW_0$ is open and non-empty.
 Let
$$
\cU_0=\cU\setminus\overline{\cW_0}.
$$
Note that $\cU_0\cup\cW_0$ is an open and dense subset of $\cU$.

%\begin{clai}\label{c.semlabel} 
%Let $\cG_0$ be  the residual subset of $\diffM$ in
%Lemma~\ref{l.robustclass}. Then for every $g\in\cG_0\cap \cU_0$,
%the class $C(P_g,g)$ does not contain any periodic point $Q_g$ with
%$\inds(Q_g)>\inds(P_g)$.
%\end{clai}
%\begin{demo} Otherwise, by definition of $\cG_0$ and Lemma~\ref{l.robustclass},
%the saddle $Q_g$ belongs robustly to $C(P_g,g)$.  Therefore, the
%diffeomorphism $g$ belongs to $\cW_0$, contradicting the
%definition of $\cU_0$. \marg{donde?}
%\end{demo}

Let $\cU_1\subset \cU$ be the set of diffeomorphisms $f$ such that
$C(P_f,f)$ has a dominated splitting $E\oplus_{_<} F$ with $\dim
(E)=\inds (P_f)$. Since the map $f\mapsto C(P_f,f)$ is
upper-semi-continuous and a dominated splitting persists in a
neighborhood of $C(P_f,f)$ by perturbations (for instance, see \cite[Chapter
B.1]{BDVbook}), one gets that the set $\cU_1$ is open. Let
$$
\cW_1=\cU \setminus \overline{\cU_1}.
$$
Then the set $\cU_1 \cup \cW_1$ is open and dense in $\cU$. As a
consequence, the open sets $\cU_0\cup \cU_1$ and $\cW_0\cap \cW_1$
are disjoint and their union is dense in $\cU$. 

Note that we are
interested in the subset of $\cU$ of diffeomorphisms $f$ whose
chain recurrence class $C(P_f,f)$ contains points of different
indices and has not an appropriate dominated splitting, that is, the set
set $\cW_0\cap \cW_1$. Thus
next lemma implies the proposition.

\begin{lemm}
\label{l.openanddense} The set $\cT$ of diffeomorphisms $g$ having a
hyperbolic set $\Si_g$ containing $P_g$  and with a $C^1$-robust
homoclinic tangency is open and dense in $\cW_0\cap\cW_1$.
\end{lemm}

\begin{demo} It is enough to prove the density of the set $\cT$. Let
$g\in\cW_0\cap\cW_1$. As $g\in\cW_0$ there is a saddle $Q_g$ with
$\inds(Q_g)>\inds(P_g)$ such that $Q_g$ belongs robustly to
$C(P_g,g)$. After a  $C^1$-perturbation, we can
assume that the diffeomorphism $g$ simultaneously belongs to the
residual set $\cG$ where $C(P_g,g)=H(P_g,g)$  and  to the residual
set $\cR$ in Theorem~\ref{th.genericblender}. Note that:

\begin{itemize}
\item
By Theorem~\ref{th.genericblender}, the set of diffeomorphisms
$g$ having a $\cut$-blender horseshoe $\La_g$ which is contained
in a transitive hyperbolic set containing $P_g$ is open and dense
in $\cW_0$.
\item
As $g\in\cW_1$ and $H(P_g,g)=C(P_g,g)$, one has that the
stable/unstable splitting defined over the set of periodic points
homoclinically related with $P_g$ is not dominated. Otherwise,
this dominated splitting could be extended to the closure of these
points (the whole $H(P_g,g)$) in a dominated way (see
\cite[Chapter B.1]{BDVbook}), which is a contradiction.
\item
Since the stable/unstable splitting defined over the set of
periodic points homoclinically related with $P_g$ is not
dominated, Theorem~\ref{t.gourmelon} implies that there is
a diffeomorphism $h$ arbitrarily $C^1$-close to $g$ with a
homoclinic tangency associated to $P_h$.
\end{itemize}
Theorem~\ref{th.zeghib} now implies that there is a
diffeomorphism $\vfi$ arbitrarily close to $g$ with a transitive
hyperbolic set containing $P_\vfi$ and having a robust homoclinic
tangency. This ends the proof of the lemma.
\end{demo}

The proof of Proposition~\ref{p.main}  is now complete.
\end{demo}

\subsubsection{End of the proof of Theorem~\ref{t.main}} 
The proof of Theorem~\ref{t.main} using Proposition~\ref{p.main} is
almost identical to the proof of Theorem~\ref{th.genericblender}
using Proposition~\ref{pr.genericblender}. It is enough to see that, for every
$n\in\NN$, there is a residual set $\cG_{\leq n}$ of
diffeomorphisms $f$ for which the conclusion of the theorem holds
for the points in $\Perfn$.

This proof is similar to the one of
Theorem~\ref{th.genericblender}. So we will omit some details.
As in Theorem~\ref{th.genericblender}, we consider the
$C^1$-open an dense subset $\cO_n\subset\diffM$ of diffeomorphisms
$f$ such that every periodic point in $\Perfn$ is hyperbolic.
Recall that the number of elements of $\Perfn$ is finite and
locally constant, and that each periodic point in $\Perfn$ has a
hyperbolic continuation on each (open) connected component of
$\cO_n$.

To state  Theorem~\ref{t.main} for periodic points in $\Perfn$,
it is enough to prove it in each connected component $\cU\in
\mbox{cc}(\cO_n)$ (recall that $\mbox{cc}(\cO_n)$ is countable):
for each connected component $\cU$, we  construct a residual
subset $\widetilde{\cG_\cU}$ such that the conclusion holds in the set
$\widetilde{\cG_\cU}\cap \cU$. Then we let
$$
\cG_\cU=\widetilde{\cG_\cU}\cup (\cO_n\setminus \cU)
$$
and define
$$\cG_{\leq n}=\bigcap_{\cU\in \mbox{cc}(\cO_n)}\cG_\cU.$$

To define $\widetilde{\cG_\cU}$ for 
a component $\cU$  of $\cO_n$, given
$f\in \cU$ write $\Perfn=\{P_{1,f}, \dots, P_{k,f}\}$
($k=k(\cU)$) and consider the continuous maps $f\mapsto P_{i,f}$,
$i\in\{1,\dots, k\}$, defined on $\cU$. For each $i$,
Proposition~\ref{p.main} provides a residual subset where the
conclusion of the theorem holds for $P_{i,f}$.  Now it is enough
to define $\widetilde{\cG_\cU}$  as the intersection of these residual sets. The proof of
Theorem~\ref{t.main} is now complete. {\hfill$\Box$}

\subsection{Robust cycles in non-hyperbolic
chain recurrence classes} \label{ss.complemet}

We close this paper by stating  an extension of \cite[Theorem
1.16]{BDijmj}. The novelty of this version is that the
hyperbolic sets involved in the robust cycle are contained in a prescribed
chain recurrence class. We note that \cite{BDijmj} does not give information about 
the relation between the
hyperbolic set involved in the robust cycle and the saddles in the
initial cycle.

\begin{theo}\label{t.heterocycle} There is a residual subset
$\cR\subset \diffM$ with the following property.  Consider any
diffeomorphism $f\in\cR$ having a chain recurrence class $C$ with
two saddles $P$ and $Q$ such that $\inds(P)=\inds (Q)+1$. Then $f$
has a $C^1$-robust heterodimensional cycle associated to
hyperbolic sets $\La$ and $\Si$ containing $P$ and $Q$.
\end{theo}

Since this result follows from arguments similar to the ones in
the previous sections and as robust heterodimensional cycles is
not the main topic of this paper, we just give some hints for the
proof.

As in the proofs above, it is enough to  state  a local version of
the theorem for a given saddle $P$. Then the general version
follows using standard genericity arguments identical to the ones in Sections~\ref{ss.generation} and \ref{s.homoclinic}.

To get  the local version of the theorem, note that for  generic diffeomorphisms $f$, the
saddle
 $Q$ is robustly in the chain recurrence class $C(P,f)$ and there is
a $\cut$-blender-horseshoe $\Si$ associated to $Q$.
In this step, the strongly intermediate points given by Proposition~\ref{p.coindexblenders} play a
key role. Finally, a
perturbation gives a robust cycle with a transverse intersection
between $W^\st(P)$ and $W^\ut(Q)$ and a robust intersection of
$W^\ut(P)$ with $W^\st(\Si)$. This completes the brief skecth of the proof.

\vskip 1cm

\flushleft{\bf Christian Bonatti}
 \ \ (bonatti@u-bourgogne.fr)\\
Institut de Math\'ematiques de Bourgogne\\ B.P. 47 870\\
21078 Dijon Cedex \\ France
\medskip

\flushleft{\bf Lorenzo J. D\'\i az}
 \ \ (lodiaz@mat.puc-rio.br)\\
Depto. Matem\'{a}tica, PUC-Rio \\ Marqu\^{e}s de S. Vicente 225\\
22453-900 Rio de Janeiro RJ \\ Brazil

\end{document}